\documentclass[twoside,11pt]{article}
\title{} \author{} \date{}
\usepackage{latexsym,amssymb,times}
\input amssym.def
\newtheorem{te}{Theorem}[section]
\newtheorem{prop}[te]{Proposition}

\newtheorem{fac}[te]{Fact}
\newtheorem{lem}[te]{Lemma}

\newtheorem{ex}[te]{Example}

\newtheorem{cla}[te]{Claim}
\def\dok{\noindent{\bf Proof. }}
\def\kdok{\hfill $\Box$ \par \vspace*{2mm} }
\def\a{\alpha}
\def\b{\beta}
\def\g{\gamma}
\def\d{\delta}

\def\o{\omega}
\def\k{\kappa}

\def\s{\sigma}

\def\l{\lambda}

\def\h{{\mathfrak h}}

\def\fh{{\mathfrak h}}

\def\fc{{\mathfrak c}}

\def\A{{\mathbb A}}
\def\B{{\mathbb B}}
\def\S{{\mathbb S}}
\def\P{{\mathbb P}}
\def\Q{{\mathbb Q}}
\def\R{{\mathbb R}}
\def\B{{\mathbb B}}
\def\N{{\mathbb N}}
\def\X{{\mathbb X}}
\def\Y{{\mathbb Y}}

\def\L{{\mathbb L}}

\def\I{{\mathcal I}}

\def\CP{{\mathcal P}}

\def\CD{{\mathcal D}}

\def\la{\langle}
\def\ra{\rangle}

\def\type{\mathop{\mbox{type}}\nolimits}
\def\Lim{\mathop{\rm Lim}\nolimits}
\def\Fin{\mathop{\rm Fin}\nolimits}
\def\sm{\mathop{\rm sm}\nolimits}
\def\sq{\mathop{\rm sq}\nolimits}
\def\asq{\mathop{\rm asq}\nolimits}
\def\rp{\mathop{\rm rp}\nolimits}
\def\supp{\mathop{\rm supp}\nolimits}
\def\red{\mathop{\rm red}\nolimits}

\def\cf{\mathop{\rm cf}\nolimits}
\def\id{\mathop{\rm id}\nolimits}
\def\cc{\mathop{\rm cc}\nolimits}
\def\Col{\mathop{\rm Col}\nolimits}

\def\Emb{\mathop{\rm Emb}\nolimits}

\def\ro{\mathop{\rm ro}\nolimits}

\def\tp{\mathop{\mathrm{tp}}\nolimits}
\def\Ord{\mathop{\mathrm{Ord}}\nolimits}

\def\Fn{\mathop{\mathrm{Fn}}\nolimits}

\def\Aut{\mathop{\mathrm{Aut}}\nolimits}
\begin{document}
\thispagestyle{plain}
\begin{center}
           {\large \bf {\uppercase{Forcing With Copies of Uncountable Ordinals}}}
\end{center}
\begin{center}
{\bf Milo\v s S.\ Kurili\'c}\footnote{Department of Mathematics and Informatics, Faculty of Sciences, University of Novi Sad,
                                             Trg Dositeja Obradovi\'ca 4, 21000 Novi Sad, Serbia.
                                             e-mail: milos@dmi.uns.ac.rs}
\end{center}
\begin{abstract}
\noindent
The set of copies ${\mathbb P} (\alpha ):=\{ f[\alpha]: f\in \mathop{\rm Emb}\nolimits (\alpha)\}$ of an ordinal $\alpha$
reconstructs the monoid $\mathop{\rm Emb}\nolimits (\alpha)$ of self embeddings of $\alpha$
and $\langle {\mathbb P} (\alpha ) ,\subset\rangle  \cong \langle \mathop{\rm Emb}\nolimits (\alpha ), \succcurlyeq ^R  \rangle$,
where $\preccurlyeq ^R$ is right Green's preorder.
Regarding the classification of countable ordinals up to isomorphism of the Boolean completions of their posets of copies
we have $\mathop{\rm ro}\nolimits (\mathop{\rm sq}\nolimits ({\mathbb P} (\alpha )))\cong \mathop{\rm ro}\nolimits ((P(\omega )/\mathop{\rm Fin})^+ \ast \pi)$,
where $\pi$ is a $P(\omega )/\mathop{\rm Fin}$-name for a $\sigma$-closed separative poset;
in particular, under ${\mathfrak h}=\omega _1 $ we have
$\mathop{\rm ro}\nolimits (\mathop{\rm sq}\nolimits ({\mathbb P} (\alpha )))\cong \mathop{\rm ro}\nolimits (P(\omega )/\mathop{\rm Fin})$,
whenever $\alpha <\o _1$.

Here we consider ordinals $\alpha \geq \omega _1$ and show that the poset $\mathop{\rm sq}\nolimits ({\mathbb P} (\alpha ))$
is either $\sigma$-closed and completely embeds a finite power of $(P(\omega )/\mathop{\rm Fin})^+$,
or completely embeds  $(P(\kappa)/[\kappa ]^{<\kappa})^+$, for some regular cardinal $\kappa >\omega$ and collapses $\omega _2$ to $\omega$.
It is consistent that  the algebra $\mathop{\rm ro}\nolimits(\mathop{\rm sq}\nolimits ({\mathbb P}(\alpha )))$
is isomorphic to $\mathop{\rm Col}\nolimits (\omega ,{\mathfrak c})$ or to $\mathop{\rm Col}\nolimits (\omega _1,{\mathfrak c})$, for each infinite ordinal $\alpha <{\mathfrak c}=\omega _2$.
More generally, the partial order $\mathop{\rm sq}\nolimits ({\mathbb P} (\alpha ))$
is isomorphic to a forcing product of the algebras of the form $P(\omega ^\delta )/{\mathcal I}_\delta$
and we detect several classes of ordinals $\delta$ such that the algebra $\mathop{\rm ro}\nolimits (P(\omega ^\delta )/{\mathcal I}_\delta)$ is isomorphic to a collapsing algebra.


\end{abstract}
\section{Introduction}\label{S1}
In this paper we continue the investigation of the posets of copies of ordinals initiated in \cite{Kord}.
Generally speaking, if $\X$ is a relational structure
and $\Emb (\X )$ is the set of its self-embeddings,
then the partial order $\la \P (\X ) ,\subset \ra$,
where $\P (\X ):=\{ f[X]: f\in \Emb (\X )\}$,
is called the {\it poset of copies of} $\X$.
Denoting by $\sq (\P )$ and $\ro(\sq (\P ))$ the separative quotient and the Boolean completion of a partial order $\P$,
the implications
$$
\X \cong \Y \Rightarrow \P(\X )\cong \P(\Y) \Rightarrow \sq \P(\X )\cong \sq \P(\Y) \Rightarrow \ro\sq \P(\X )\cong \ro\sq \P(\Y)
$$
express the relationship between ``similarity relations" in the class of relational structures
and each of them provides a classification of structures \cite{Ktow,Kdif}.
The last one, $\ro\sq \P(\X )\cong \ro\sq \P(\Y)$,
gives the coarsest classification,
it is equivalent to the forcing equivalence of posets of copies,
and, hence, it can be explored using forcing-theoretic methods.
The related results can be regarded as the statements belonging to the theory of relations, set theory or the theory of Boolean algebras.

In addition, these results can be applied in investigation of the transformation monoids
$\la \Emb (\X ), \circ , \id _X\ra$.
First we note that
the set $\Aut (\X)$ of automorphisms of $\X$
together with a choice of an embedding
$f_A \in \Emb (\X )$ such that $f_A [X]=A$, for each $A\in \P (\X )$,
determine the set of embeddings: $\Emb (\X ) = \{ f_A \circ a : A\in \P (\X ) \mbox{ and } a\in \Aut (\X )\}$.
Second, we recall that {\it right Green's preorder} on $\Emb (\X )$ is defined by
$f \preccurlyeq ^R g $ iff $f\circ h =g $, for some $h \in \Emb (\X )$,
and note that, defining $f \thickapprox ^R g$ iff $f \preccurlyeq ^R g$ and $g \preccurlyeq ^R f$, we have
$f \preccurlyeq ^R g  \Leftrightarrow  g[X]\subset f[X]$ and $f \thickapprox ^R g\Leftrightarrow  g[X]= f[X]$.
So, $\la \P (\X ), \subset\ra$ is the antisymmetric quotient of the preorder $\la \Emb (\X ) ,\succcurlyeq ^R \ra$.
If, in particular, $\X$ is an ordinal $\a $,
it is a rigid structure and for each copy $A\in \P (\a )$
there is a {\it unique} $f_A \in \Emb (\a )$ such that $f_A [\a]=A$.
Thus, here we have $\Emb (\a ) = \{ f_A : A\in \P (\a ) \}$
and defining the binary operation $\ast$ on $\P (\a )$ by $A\ast B =f_A[B]$
and the mapping $F:\P (\a ) \rightarrow \Emb (\a )$ by $F(A)=f_A$,
we obtain an isomorphism of monoids; namely,
$F:\la \P (\a ),\ast ,\a \ra \cong \la \Emb (\a ), \circ , \id _\a\ra $.
In addition, the preorder  $\preccurlyeq ^R$ on $\Emb (\a )$ becomes a partial order and the same mapping $F$ witnesses that
$\la \P (\a ) ,\subset\ra  \cong \la \Emb (\a ), \succcurlyeq ^R \ra $.

It turns out that, regarding countable linear orders,
the main roles in the description of the algebra $\ro\sq \P(\X )$
are played by the poset $\S$ of perfect subsets of the real line ordered by inclusion (the Sacks forcing)
and the algebra $P(\o )/\Fin$.
Namely, if $\X$ is non-scattered (i.e.,  $\Q \hookrightarrow \X$),
then $\ro\sq \P(\X )\cong \ro\sq (\S \ast \pi)$, where $\pi$ is an $\S$-name for a $\s$-closed poset and, consistently, the name for $P(\o )/\Fin$ in the Sacks extension \cite{KurTod}.
On the other hand, if  $\X$ is a scattered linear order, then the poset $\sq \P(\X )$ is $\s$-closed, atomless and, consistently, $\ro\sq \P(\X )\cong \ro\sq (P(\o )/\Fin)$ \cite{Kscatt}.
We note that for the real line $\R$ we have $\ro\sq \P(\R )\cong \ro\sq (\S )$ \cite{Kmono}
and that the Sacks forcing appears as a factor in the presentation of posets of copies of several ultrahomogeneous structures,
 \cite{GK,KurTod,KTR1,KTR,KTlab,KTT}.

In \cite{Kord} the posets of copies of countable ordinals $\a \geq \o$ were investigated
and, in order to present the main results, we explain some notation.
If $\B$ is Boolean lattice, $\rp (\B )$ denotes the reduced power $\B ^\o / \Fin$,
and the iterated reduced powers $\rp ^n(\B )$, for $n\in \o$, are defined by $\rp ^0(\B) =\B $ and $\rp ^{n+1}(\B )= \rp (\rp ^n (\B ))$.
If $\d >0$ is an ordinal, then $\I_{\o ^\d}:=\{ I\subset \o ^\d :\o ^\d \not\hookrightarrow I \}\subset P(\o ^\d)$ is an ideal
and the Boolean algebra $P(\o ^\d)/ \I_{\o ^\d}$ is the corresponding quotient.
The first result shows that the separative quotient of the poset $\P (\alpha )$
is isomorphic to the forcing product of finitely many Boolean algebras:
if $\a=\o ^{\g _n +r_n }s_n + \dots + \o ^{ \g _0 +r_0 }s_0 +k$ is a countable ordinal presented in the Cantor normal form,
where $k\in \o$, $r_i \in \o$, $s_i \in \N$, $\g _i \in \Lim \cup \{ 1 \}$
and $\g _n +r_n > \dots > \g _0 +r_0$, then
\begin{equation}\label{EQ253}\textstyle
\sq (\P (\a ))\cong \prod _{i=0}^n \Big( \Big( \rp ^{r_i}( P(\o ^{\g _i} )/ \I _{\o ^{\g _i} })\Big)^+ \Big)^{s_i} .
\end{equation}
The second result is related to the aforementioned coarse classification of ordinals
and their embedding monoids.
We recall that the {\it distributivity number} $\fh (\P )$ of a partial order $\P$
is the minimal size of a family $\CD$ of open dense subsets of $\P$ such that $\bigcap \CD$ is not dense,
and that $\o_1 \leq \fh :=\fh ((P(\o )/\Fin )^+)\leq \fc$.
\begin{te}\label{T4421}
If $\a $ is a countable ordinal, then

{\rm (a)} $\ro (\sq (\P (\a )))\cong \ro ((P(\o )/\Fin )^+ \ast \pi)$, where $[\o ] \Vdash$ ``$\pi$ is a $\s$-closed separative poset";

{\rm (b)} $\ro (\sq (\P (\a )))\cong \ro ((P(\o )/\Fin )^+)$, if $\,{\mathfrak h}=\o _1$.
\end{te}
Thus if ${\mathfrak h}=\o _1$, then all countable ordinals belong to the same class with respect to the coarse classification.
We note that the equality ${\mathfrak h}=\o _1$ follows from CH
and holds in many models of $\neg\,$CH;
for example, in all iterated forcing models obtained by adding Cohen, random, Sacks, Miller, Laver and Hechler reals,
(see \cite{Blas}).
But the statement that $\ro (\sq (\P (\a )))\cong \ro ((P(\o )/\Fin )^+)$ holds for all countable ordinals $\a$
is independent of ZFC (see \cite{Kord}).

In this paper we investigate the posets $\P (\a )$ for uncountable ordinals $\a$.
Section \ref{S2} contains definitions and known facts used in the paper.
In Section \ref{S3} we show that the separative quotient $\sq \P (\a )$
is isomorphic to the direct product of $\sum _{i=1}^n s_i$-many factors $\sq \P (\o ^{\d _i})$,
where $\a = \o ^{\d _n}s_n+\dots +\o ^{\d _1}s_1+ m$ is the Cantor normal form for $\a$,
and reduce the analysis to the investigation of these factors.

So, in Section \ref{S4} we regard the posets of the form $\P (\o ^\d )$.
It turns out that the behavior of $\P (\o ^\d )$ depends on the cofinality $\k :=\cf (\d )$.
If $\k \leq \o$ (Case (A)), then the poset $\sq (\P (\o ^\d ))$ is $\s$-closed, completely embeds $(P(\o )/\Fin )^+$ and, hence, collapses $\fc$ to $\fh$.
Otherwise we have  $\k \geq \o _1$ and, presenting $\d$ in its Cantor normal form in the base $\k$,
we detect additional four cases, (B)--(E).
If (B) holds, then again, $\sq (\P (\o ^\d ))$ is $\s$-closed and completely embeds $(P(\o )/\Fin )^+$.
In Cases (D) and (E) the poset $\sq (\P (\o ^\d ))$ completely embeds the poset $\CP _\k :=(P(\k )/[\k ]^{<\k })^+ $
and collapses (at least) all cardinals $\mu < \cc (\CP _\k)$ (and, hence, $\o _2$) to $\o$.
Finally, in Case (C)  $\sq (\P (\o ^\d ))$ completely embeds $\CP _\l $, where $\o <\l <\k$,
and collapses (at least) all cardinals $\mu < \cc (\CP _\l)$ (and, again, $\o _2$) to $\o$.

Thus, by the (ZFC) results of Section \ref{S4}, the poset $\P (\o ^\d )$ either preserves $\o _1$ and forces $|\fc|=|\fh|$,  or collapses (at least) $\o _2$ to $\o$.
So, it is natural to ask whether $\P (\o ^\d )$ collapses more cardinals
and a more ambitious task is to characterize the algebra $\ro (\sq (\P (\o ^\d )))$.
Since $|\sq (\P (\o ^\d ))|=|P(\o ^{\d })/ \I_{\o ^{\d }}| \leq 2^{|\d|}$, the extreme case is when
$\P (\o ^\d )$ collapses $2^{|\d|}$ to $\o _1$, in cases (A) and (B),
or collapses $2^{|\d|}$ to $\o $, in other cases.
So, in Section \ref{S5} we consider that phenomenon of ``maximal collapse".
For example, in cases (D) and (E) we have $\ro (\sq (\P (\o ^\d )))\cong \Col (\o ,2^{|\d|} )$, if $2^{\cf (\d )} =2^{|\d|}$ and
($2^{<\cf (\d )}=\cf (\d )$ or $2^{\cf (\d )}=\cf (\d )^+$). Regarding case (A) we have
$\ro (\sq (\P (\o ^{\d +n})))\cong \Col (\o _1, 2^{|\d |})$, for all $n\in \N$, whenever $\P (\o ^\d)$ collapses $2^{|\d |}$ to $\o$,
or $\sq(\P (\o ^\d))$ is $\s$-closed and collapses $2^{|\d |}$ to $\o_1$.
\section{Preliminaries }\label{S2}
\paragraph{Partial orders and forcing}
If $\P =\la P, \leq \ra $ is a preorder, the elements $p$ and $q$ of $P$ are {\it incompatible},
we write $p\perp q$, iff there is no $r\in P$ such that $r\leq p,q$.
A set $A\subset P$ is an antichain iff $p\perp q$, for different $p,q\in A$.
An element $p$ of $P$ is an {\it atom } iff each two elements $q,r\leq p$ are compatible.
${\mathbb P} $ is called: {\it atomless} iff it has no atoms;
{\it homogeneous} iff  it has a largest element and $\P \cong (\cdot ,p] $, for each $p\in \P$.
A set $D\subset P$ is called {\it dense} iff for each $p\in P$ there is $q\in D$ such that $q\leq p$;
$D$ is called {\it open} iff $q\leq p\in D$ implies $q\in D$.
If $\kappa $ is a cardinal, ${\mathbb P} $ is called $\kappa ${\it -closed} iff for each
$\gamma <\kappa $ each sequence $\langle  p_\alpha :\alpha <\gamma\rangle $ in $P$, such that $\alpha <\beta \Rightarrow p_{\beta}\leq p_\alpha $,
has a lower bound in $P$. $\omega _1$-closed preorders are called {\it $\sigma$-closed}.

A preorder ${\mathbb P} =\langle  P , \leq \rangle $ is called
{\it separative} iff for each $p,q\in P$ satisfying $p\not\leq q$ there is $r\leq p$ such that $r \perp q$.
The {\it separative modification} of ${\mathbb P}$
is the separative preorder $\mathop{\rm sm}\nolimits ({\mathbb P} )=\langle  P , \leq ^*\rangle $, where
$p\leq ^* q \Leftrightarrow \forall r\leq p \; \exists s \leq r \; s\leq q $.
The {\it separative quotient} of ${\mathbb P}$
is the separative partial order $\mathop{\rm sq}\nolimits  ({\mathbb P} )=\asq (\sm (\P))$, where $\asq (\P )$ denotes the {\it antisymmetric quotient} of a preorder $\P$.
Preorders $\P$ and $\Q$ are called {\it forcing equivalent}, in notation $\P \equiv_{forc}\Q$, iff they produce the same generic extensions.
\begin{fac}(Folklore)  \label{T2226}
If ${\mathbb P} $ and ${\mathbb P} _i$, for $i\in I$, are preorders, then

(a) ${\mathbb P}\equiv_{forc}\mathop{\rm sm}\nolimits ({\mathbb P})\equiv_{forc}\mathop{\rm sq}\nolimits  ({\mathbb P})\equiv_{forc}\ro(\mathop{\rm sq}\nolimits  ({\mathbb P}))$;

(b) If ${\mathbb P} _i$, $i\in I$, are $\kappa $-closed preorders, then $\prod _{i\in I}{\mathbb P} _i$ is $\kappa $-closed;

(c) $\mathop{\rm sm}\nolimits (\prod _{i\in I}{\mathbb P} _i) = \prod _{i\in I}\mathop{\rm sm}\nolimits {\mathbb P} _i$ and
 $\mathop{\rm sq}\nolimits  (\prod _{i\in I}{\mathbb P} _i) \cong \prod _{i\in I}\mathop{\rm sq}\nolimits  {\mathbb P} _i$;

(d) If $X$ is an infinite set, $\I \subset P(X)$ is an ideal containing $[X]^{<\o }$,
$\I ^+ = P(X)\setminus \I$ is the corresponding co-ideal and $A\subset _\I B$ iff $A\setminus B \in \I$, for $A,B \in \I ^+$, then
$\sm \la \I ^+ , \subset \ra $ $= \la \I ^+ , \subset _\I \ra$ and $\sq \la \I ^+ , \subset \ra =(P(X)/\I )^+$.
\end{fac}
\noindent
If $\la \P ,\leq _\P , 1_\P \ra$ and  $\la \Q , \leq _\Q , 1_\Q \ra$ are preorders, then a mapping
$f: \P \rightarrow \Q$ is a {\it complete embedding}, in notation $f:\P \hookrightarrow _c \Q$  iff

(ce1) $p_1 \leq _\P p_2 \Rightarrow f(p _1)\leq _\Q f(p_2)$,

(ce2) $p_1 \perp _\P p_2 \Leftrightarrow f(p _1)\perp _\Q f(p_2)$,

(ce3) $\forall q\in \Q \; \exists p\in \P \; \forall p' \leq _\P p \;\; f(p') \not\perp _\Q q$.

\noindent
Then, for $q\in \Q $,
$\red _\P (q) =\{ p\in \P : \forall p' \leq _\P p \;\; \exists q'\in \Q \;\;q' \leq f(p') , q \}$
is the set of {\it reductions} of $q$ to $\P$. The following fact is folklore (see \cite{Kun}).
\begin{fac}\label{T4432}
If $f:\P \hookrightarrow _c \Q$, then $\Q$ is forcing equivalent to the two step iteration $\P \ast \pi $, where
$\pi $ is a $\P$-name for a preorder.
In addition, $F:\asq (\P) \hookrightarrow _c \asq (\Q )$, where $F([p]_{\asq (\P)})=[f(p)]_{\asq (\Q)}$, for $p\in \P$.
\end{fac}
\begin{fac}\label{T227}
If ${\mathbb P} _i \hookrightarrow _c{\mathbb P} _i'$, for $i\in I$, then $\prod _{i\in I}{\mathbb P} _i\hookrightarrow _c \prod _{i\in I}{\mathbb P} _i'$.
\end{fac}
\begin{fac}\label{T225}
If $V\models \mathrm{GCH}$
and, in $V$, $\k>\o$ is a regular cardinal and $\P =\Fn (\k ,2)$ is the Cohen forcing,
then in $V[G]$ we have $\fc =\k$.
If, in $V$, $\mu$ and $\theta$ are infinite cardinals, then
$(\theta ^\mu)^{V[G]}=(\max(\k ,\theta)^\mu)^V $ (see \cite{Kun}, p.\ 246).
\end{fac}
\paragraph{Collapsing algebras and the algebras $P(\k )/[\k ]^{<\k}$}
Let $\k\geq 2$ and $\l \geq \o$ be cardinals.
The {\it collapsing algebra} $\Col (\l ,\k)$ is the Boolean completion  of the reversed tree $\la {}^{<\l }\k ,\supset\ra$.
The following fact is folklore; see \cite{KCol}, Thm 3.1 for a proof.
\begin{fac}\label{T4130}
 Let $\l\geq \o$ be a regular cardinal and $\P$ a separative $\l$-closed preorder of size $\k=\k^{<\l}$.

(a) If $\k >\l$ and $1_\P\Vdash |\check{\k}|= \check{\l}$, then $\ro (\sq (\P ))\cong \Col (\l ,\k)$;

(b) If $\k =\l$ and $\P$ is atomless, then $\ro (\sq (\P ))\cong \Col (\l ,\l)$.
\end{fac}
For a cardinal $\k\geq \o$ the partial order $(P(\k )/[\k ]^{<\k })^+$ is denoted by $\CP _\k $;
thus, $\CP _\o =(P(\o )/\Fin )^+ $.
If $\subset ^*$ is the preorder on the set $[\k ]^\k$ given by $A\subset ^* B$ iff $|A\setminus B|<\k$,
then, clearly, $\sm \la [\k ]^\k ,\subset \ra = \la [\k ]^\k ,\subset ^*\ra$ and $\sq \la [\k ]^\k ,\subset \ra = \CP _\k $.
\begin{fac}\label{T209}
If $\k$ is an infinite cardinal, then for the partial order $\CP _\k$ we have

(a) $\CP _\k$ is atomless, homogeneous, of size $2^\k$ and $\k ^{++}\leq \cc (\CP _\k)\leq (2^\k) ^+$;

(b)  Forcing by $\CP _\o$ collapses $\fc$ to $\h$ (Balcar, Pelant and Simon, \cite{BPS});

(c)  If $\k >\cf (\k ) =\o$, then $\Col(\o _1, \k ^\o)\hookrightarrow _c \ro (\CP _\k)$ (Kojman, Shelah \cite{Koj});

(d)  If $\k > 2^{\cf (\k )}> \cf (\k ) >\o$, then  $\Col(\o , \k ^+)\hookrightarrow _c \ro (\CP _\k)$ (Shelah, \cite{She});

(e)  If $\k = \cf (\k )>\o$, then $\CP _\k $ collapses each $\lambda <\cc (\CP _\k)$ to $\o$ (Shelah, \cite{She}).
\end{fac}
\paragraph{Linear orders. Indecomposable ordinals}
A linear order $\la L, < \ra$ is called {\it dense} iff for each $x,y\in L$ satisfying $x<y$ there is $z\in L$ such that
$x<z<y$. $\la L, < \ra$ is said to be {\it scattered} iff it does not contain a dense suborder iff $\Q \not\hookrightarrow L$.
\begin{fac}  \label{T4102}
A linear order $L$ satisfying $L+L\hookrightarrow L$ is not scattered (\cite{Rosen}, p.\ 176).
\end{fac}
\begin{fac}\label{T4400}
For a limit ordinal $\a$ the following conditions are equivalent:

(a) $\b + \g < \a$, for each $\b , \g <\a$, (that is, $\a $ is  indecomposable),

(b) $\b +\g= \a  \land \g>0 \Rightarrow \g =\a$ (that is, $\a $ is right indecomposable),

(c) $\,\b +\a =\a$, for each $\b <\a$,

(d) $\a = \o ^\delta$, for some ordinal $\delta >0$,

(e) $\a$ is an indivisible structure ($\a\hookrightarrow A$ or $\a\hookrightarrow B$, whenever $\a =A\;\dot{\cup}\; B$),

(f) $\I _\a =\{ I\subset \a : \a \not\hookrightarrow I\}$ is an ideal in $P(\a )$.
\end{fac}
\dok
For (a) $\Leftrightarrow $ (d) $\Leftrightarrow $ (c) see \cite{Sier}, p.\ 282, 323. For (b) $\Leftrightarrow $ (d) see \cite{Rosen}, p.\ 176.
The equivalence (a) $\Leftrightarrow $ (e) is 6.8.1 of \cite{Fra}. (e) $\Leftrightarrow$ (f) is evident.
\hfill $\Box$
\begin{fac}\label{T4128}
(a) $\o ^\delta \geq \delta$, for each ordinal $\d$;

(b) $|\o ^\delta |=\o |\delta|$, for each $\d \geq 1$;

(c) If $\a = \o ^{\d _n}s_n+\dots +\o ^{\d _1}s_1+ m$ is an uncountable ordinal, then $|\a |=|\d _n|$.
\end{fac}
\dok
(a) First we have $\o ^0 =1 \geq 0$.
Suppose that $\o ^\delta \geq \delta$, for each ordinal $\d <\theta$.
If $\theta =\d +1$, then $\o ^\theta =\o ^\d \o >\o ^\d 1=\o ^\d \geq \d$ and, hence, $\o ^\theta \geq\d +1=\theta$.
If $\theta$ is a limit ordinal, then $\o ^\theta =\sup \{ \o ^\d :\d <\theta \}\geq \sup \{ \d :\d <\theta \}=\theta$.

(b) If $\d \geq 1$, then $|\o ^\delta |\geq \o$ and, by (a), $|\o ^\delta| \geq |\delta|$, so $|\o ^\delta |\geq \o |\delta|$. Let
\begin{equation}\label{EQ4163}
\forall \d <\theta \;\;|\o ^\delta |=\o |\delta| .
\end{equation}
If $\theta =\d +1$, then $|\o ^\theta |=|\o ^\d \o |=|\o ^\d|\o=|\delta|\o\o =|\delta +1|\o=|\theta |\o$.

If $\theta$ is a limit ordinal, then $\o ^\theta =\sup \{ \o ^\d :\d <\theta \}$
and by (\ref{EQ4163}), for each $\d <\theta$ we have $|\o ^\delta |=\o |\delta|\leq \o |\theta |=|\theta |$.
So the ordinal $\o ^\theta $ is an union of $|\theta |$-many ordinals of size $\leq |\theta |$
and, hence, $|\o ^\theta |=|\theta|=\o|\theta|$.

(c) $\a$ is a finite sum of summands of the form $\o ^{\d _i}$ and by (b) we have $|\o ^{\delta _i}|=\o |\delta _i|$.
So $\a \geq \o _1$ implies that the largest summand, $\o ^{\d _n}$ is uncountable and, hence, $|\d _n|\geq \o _1$.
Since $i<n$ implies $\d _i <\d _n$ we have
$\o ^{\d _i}<\o ^{\d _n}$, so $|\o ^{\d _i}| \leq |\o ^{\d _n}|$. Thus by (b) we have $|\a|=|\o ^{\d _n}|=\o |\delta _n|=|\delta _n|$.
\hfill $\Box$
\section{Factorization of $\P (\a)$}\label{S3}
Here we start with analysis of the posets of the form $\P (\a)$, where $\a$ is an uncountable ordinal.
Clearly, these posets are homogeneous
and atomless.
\begin{lem}\label{T4409}
$\la \P (\g +m ),\subset \ra \cong \la  \P (\g ),\subset \ra$, for each limit ordinal $\g$ and each $m\in \N$.
\end{lem}
\dok
First we prove $\P (\g +m )= \{ C \cup \{ \g , \g +1, \dots , \g +m-1\} : C\in \P (\g ) \} $.
The inclusion ``$\supset$" is evident.
If $A\in \P (\g +m )$ and $f:\g +m \hookrightarrow \g +m$, where $A=f[\g +m]$,
then, since $f$ is an increasing function, we have $f(\b )\geq \b$, for each $\b \in \g +m$,
which implies $f(\g +i) =\g +i$, for $i<m$,
and, hence, $C= f[\g ] \in \P (\g )$ and $A=C \cup \{ \g , \g +1, \dots , \g +m-1\}$.

Now it is easy to show that the mapping
$F: \la \P (\g ),\subset \ra \rightarrow \la \P (\g +m ),\subset\ra$,
given by $F(C)=C \cup \{ \g , \g +1, \dots , \g +m-1\} $, is an isomorphism.
\hfill $\Box$
\begin{te}\label{T200}
If $\a = \o ^{\d _n}s_n+\dots +\o ^{\d _1}s_1+ m$ is an infinite ordinal presented in the Cantor normal form,
where $n, s_1, \dots , s_n\in \N$,  $\;0< \d_1 < \dots <\d _n$ and $\,m \in \o$, then
$\P (\a )\cong \prod _{i=1}^n \big(\P (\o ^{\d _i})\big)^{s_i}$,
$\sm \P (\a )\cong \prod _{i=1}^n \big(\sm \P (\o ^{\d _i})\big)^{s_i}$ and
$\sq \P (\a )\cong \prod _{i=1}^n \big(\sq \P (\o ^{\d _i})\big)^{s_i}$.
\end{te}
\dok
Since $\d _1,s_1 >0$, by Lemma \ref{T4409} we have $\P (\a )\cong \P (\g )$,
where, defining $k:= \sum _{i=1}^n s_i$,
we have $\g = \o ^{\d _n}s_n+\dots +\o ^{\d _1}s_1 =\sum _{j=1}^k \L _j$.
By induction on $k\in \N$ we prove that $\P (\g )\cong \prod _{i=1}^n \big(\P (\o ^{\d _i})\big)^{s_i}(\cong\prod _{j=1}^k \P (\L _j)$), which will follow from
\begin{equation}\label{EQ201}\textstyle
\P (\g )=\big\{ \bigcup _{j=1}^k C_j : \la C_1,\dots ,C_k\ra \in \prod _{j=1}^k \P (\L _j)\big\}.
\end{equation}
For $k=1$ we have $n=s_1 =1$, so $\g =\o ^{\d _1}=\L _1$, $\P(\g )=\P(\L _1)$ and (\ref{EQ201}) is true.

Let (\ref{EQ201}) be true for $k$,
let $n, s_1, \dots , s_n\in \N$, where $\sum _{i=1}^n s_i=k+1$,
let $\;0< \d_1 < \dots <\d _n$ and
$\g '= \o ^{\d _n}s_n+\dots +\o ^{\d _1}s_1 =\sum _{j=1}^{k+1} \L _j$.
Then $\g '=\L _1 +\sum _{j=2}^{k+1} \L _j$, where $\L _1\cong\o ^{\d _n}$
and (\ref{EQ201}) is true for $\g :=\sum _{i=2}^{k+1} \L _j$.

Let $C\in \P (\g ')$ and $f\in \Emb (\g ')$, where $f[\g ']=C$.
We show that $C_1:=f[L_1]$ is a cofinal subset of $\L _1$;
then we will have $f[\g ]\subset \g$ and, hence, $f[\g ]\in \P (\g )$;
so, by the induction hypothesis (\ref{EQ201}) will be true for $k+1$.

Assuming that $F:=f[L_1]\cap \g \neq \emptyset$
the order $F$ would be isomorphic to a final part of $\L _1\cong\o ^{\d _n}$
and, by Fact \ref{T4400}(b) we would have $F\cong\o ^{\d _n}$.
So we would have $\g ' \cong F +f[\g]\subset \g \subset \g '$
and, hence $\g ' \cong \g =\o ^{\d _n}(s_n-1)+\dots +\o ^{\d _1}s_1$,
which is impossible since each ordinal has a unique representation in the Cantor normal form.
Thus $f[L_1]\subset L_1$ and $f[L_1]\cong L_1\cong \o ^{\d _n}$,
which implies that $C_1:=f[L_1]$ is a cofinal subset of $\L _1$ and $C_1\in \P (\L _1)$.
Consequently, $f[\g ]\subset \g$ and, hence, $f[\g ]\in \P (\g )$;
so, by the induction hypothesis $f[\g]=\bigcup _{j=2}^{k+1} C_j$, where $\la C_2,\dots ,C_k\ra \in \prod _{j=2}^{k+1} \P (\L _j)$
and $f[\g ']=\bigcup _{j=1}^{k+1} C_j$.

Thus, $\P (\a )\cong \prod _{i=1}^n (\P (\o ^{\d _i}))^{s_i}$ and the rest follows from Fact \ref{T2226}(c).
\hfill $\Box$
\section{Factorization of factors $\P (\o ^\delta)$. ZFC results}\label{S4}
\paragraph{Classification of exponents $\delta$}
If $\d$ is an ordinal, $\cf (\d )=\k\geq \o$,
$\la \d _\xi : \xi\in \k \ra \in {}^{\k}\d $ is a strictly increasing sequence of ordinals
and $\d =\sup \{ \d _\xi : \xi\in \k\}$,
then we will write $\d=\lim _{\xi \rightarrow \k}\d _\xi$. Let $\Lim _0:=\{ 0 \}\cup \Lim$.
\begin{lem}\label{T4121}
If $\d>0$ is an ordinal and $\cf (\d )=\k $, then exactly one of the following conditions is satisfied.
\begin{itemize}
\item[(A)] $\d$ is a successor ordinal or $\cf (\d )=\o$ (that is, $\k \leq \o$);
\item[(B)] $\d =\theta +\k$, where $\Ord \ni \theta\geq \k>\cf (\theta )=\o$,\\
           and $\theta =\lim _{n\rightarrow \o}\d _n$, where $\cf (\d _n)=\k$, for all $n\in \o$;
\item[(C)] $\d =\theta +\k$, where $\Ord \ni \theta\geq \k >\cf (\theta )>\o$,\\
           and $\theta =\lim _{\xi\rightarrow \,\lambda }\d _\xi $, where $\,\cf (\d _\xi)=\k$, for all $\xi\in \lambda :=\cf (\theta )$;
\item[(D)] $\d =\theta +\k$, where  $\Ord \ni \theta\geq \cf (\theta )\geq \k >\o$ or $\theta =0$;
\item[(E)] $\d=\lim _{\xi \rightarrow \k}\d _\xi$, where $\cf (\d _\xi)=\k >\o$, for all $\xi\in \k$.
\end{itemize}
\end{lem}
\dok
Let $\k>\o$. Since $\d\geq \cf (\d )= \k$, the ordinal $\d$ can be presented in the Cantor normal form in the base $\k$
\begin{equation}\label{EQ4146}
\d = \k ^{\xi _n}\zeta _n + \k ^{\xi _{n-1}}\zeta _{n-1} + \dots +\k ^{\xi _1}\zeta _1 +\b ,
\end{equation}
where
$n \in \N$,
$1\leq \zeta _k < \k$, for $k\leq n$,
$0\leq \b <\k$ and
$1\leq \xi _1 < \xi _2 < \dots  < \xi _n $.
First, $\b >0$ would mean that $\b$, and, hence, $\d$, is a successor ordinal or a limit ordinal of cofinality $<\k$; so $\b=0$.

Second, since $\zeta _1 <\k$, $\zeta _1 \in \Lim$ would imply that $\cf (\zeta _1)<\k$
and the final part $\k ^{\xi _1}\zeta _1$ of $\d$,
which is a $\zeta _1$-sum of copies of $\k ^{\xi _1}$
would be of cofinality $<\k$;
this would imply that $\cf (\d )<\k$.
Thus $\k ^{\xi _1}\zeta _1= \k ^{\xi _1}(\g  + m )=\k ^{\xi _1}\g   + \k ^{\xi _1}+\dots +\k ^{\xi _1}$,
where $\g \in \Lim_0 \;\cap \;\k $
and $m\in \N$,
so $\d$ has a final part isomorphic to $\k ^{\xi _1}$.

Third, $\o\leq\cf (\xi _1)<\k$ would imply that $\cf (\k ^{\xi _1})<\k$ and, hence, $\cf (\d )<\k$.
So either $\cf (\xi _1)=\k$ or $\xi _1=\g ' +m'$,
where $\g '\in \Lim _0$ and $m'\in \N$.
Thus defining
\begin{equation}\label{EQ4146'}
\theta = \k ^{\xi _n}\zeta _n +  \dots + \k ^{\xi _{2}}\zeta _{2} +\k ^{\xi _1}\g   + \k ^{\xi _1}(m-1),
\end{equation}
\begin{equation}\label{EQ4146''}
\mbox{where }\;\zeta _1=\g  +m \;\;\mbox{ and }\;\; \g \in \Lim _0 \;\cap \;\k \;\;\mbox{ and }\;\; m\in \N,
\end{equation}
we have
\begin{eqnarray}
\d &=& \theta + \k ^{\xi _1}, \mbox{ where } \cf (\xi _1)=\k, \mbox{ or }   \label{EQ4147}     \\
\d &=& \theta +\k ^{\g ' +m'}, \mbox{ where }\g ' \in \Lim _0 \mbox{ and } m'\in \N.\label{EQ4148}
\end{eqnarray}
If (\ref{EQ4147}) holds,
let $\la \a _\xi : \xi <\k \ra$ be an increasing  sequence cofinal in $\xi _1$
and let us define $\d _\xi :=\theta +\k^{\a _\xi +1}$, for $\xi <\k$.
Then for $\xi <\xi '<\k$ we have $\d _\xi < \d _{\xi '}<\d$
and  $\sup \{ \d _\xi : \xi < \k \}= \theta + \sup \{ \k ^{\a _\xi +1} : \xi < \k \}=\theta + \k ^{\xi _1}=\d$.
So the sequence $\la \d _\xi : \xi\in \k \ra$ is increasing, cofinal in $\d$
and $\cf (\d _\xi)=\cf (\k^{\a _\xi}\k )=\k$,
because $\k ^{\a _\xi}\k$ is a $\k$-sum. So (E) is true.

If (\ref{EQ4148}) holds and $\g ' +m' \geq 2$,
then we have $\a:= \g ' +m' -1\geq 1$ and $\d =\theta + \k ^\a\k$.
Thus the final part $\k ^\a\k$ of $\d$ is a $\k$-sum of copies of the ordinal $\k ^\a\geq \k$.
Now defining the ordinals $\d _\xi := \theta +\k ^\a \xi + \k$, for $\xi <\k$, we have
$$
\k ^\a \xi + \k
\leq \k ^\a \xi + \k ^\a
=\k ^\a (\xi +1)
<\k ^\a (\xi +1) +\k
\leq\k ^\a (\xi +2)
<\k^\a\k
$$
and, hence, $\d _\xi<\d$.
Clearly we have $\cf (\d _\xi)=\k$, for all $\xi < \k $,
and $\la \d _\xi : \xi < \k \ra$ is  an increasing  sequence cofinal in $\d$.
So (E) is true again.

If (\ref{EQ4148}) holds and $\g ' +m' <2$,
then $\g '=0$, $m'=1$, $\xi _1=1$ and $\d =\theta +\k$,
where by (\ref{EQ4146'}) and (\ref{EQ4146''}) we have
\begin{equation}\label{EQ4159}
\theta = (\k ^{\xi _n}\zeta _n +  \dots +\k ^{\xi _2}\zeta _2) + \k(\g  +m-1) =\theta ' + \k\a,
\end{equation}
where $\theta '=0$, if $n=1$, and $\a:=\g  +m-1 $ is an ordinal satisfying $0\leq \a <\k$. If $\theta =0$, then we have (D).

If $\theta >0$, then $\a >0$ or $\theta ' >0$ (and, hence, $\xi_2 >1$ and $\zeta _2 \geq 1$), which gives $\theta \geq \k$.
In addition, the equality (\ref{EQ4159}) presents the ordinal $\theta$ in the Cantor normal form and, because of its uniqueness,
$\theta$ can not be presented in the form $\theta _1 + \b$, where $0<\b <\k$;
thus, $\theta$ is a limit ordinal.
First, if $\cf (\theta )\geq \k$, we have (D).

Second, if $\lambda :=\cf (\theta )<\k$,
let $\la \a_\xi : \xi<\lambda \ra$ be an increasing sequence cofinal in $\theta$.
Assuming that $\a _\xi + \k \geq \theta$, for some $\xi <\lambda$,
since $\cf (\theta )<\k =\cf (\a _\xi + \k)$ we would have $\a _\xi + \k > \theta$,
and, hence, $\theta =\a _\xi + \b$, where $0<\b <\k$, which is false.
So $\a _\xi + \k < \theta$, for $\xi <\lambda$,
and $\theta =\sup \{\a_\xi : \xi <\lambda \}\leq \sup \{\a_\xi + \k : \xi <\lambda \}\leq \theta$.
Defining $\d_\xi =\a_\xi + \k$, for $\xi <\lambda$,
we have $\theta =\sup \{\d_\xi : \xi <\lambda \}$ and $\cf (\d _\xi )=\k$, for all $\xi <\lambda$.
Now, if $\lambda =\o$ we have (B) and if $\lambda >\o$ we have (C).

Finally we show that each $\d$ satisfies exactly one of conditions (A)--(E).
If (A) is true, then $\k \leq \o$ and (B)--(E) are false, because (B)--(E) imply that $\k>\o$.
If (E) holds, then (B)--(D) are false, since (B)--(D) imply that $\d =\theta +\k$ and for an ordinal $\b \in (\theta, \theta +\k)$ we have $\cf (\b )<\k$.

Suppose that $\d$ satisfies (B) and (C). Then $\d =I+F=I'+F'$, where $F\cong F'\cong \k$,
$I\cong \theta\geq \k>\cf (\theta )=\o$,
and $\theta =\lim _{n\rightarrow \o}\d _n$, where $\cf (\d _n)=\k$, for all $n\in \o$;
and $I' \cong \theta '\geq \k >\cf (\theta ' )=\lambda >\o$,
and $\theta ' =\lim _{\xi\rightarrow \,\lambda }\d _\xi $, where $\,\cf (\d _\xi)=\k$, for all $\xi\in \lambda $.
Since $\cf (\theta )<\cf (\theta ')$ we have $I\neq I'$; so either $I\varsubsetneq I'$ or $I '\varsubsetneq I$.
If $I\varsubsetneq I'$, then there is $\xi\in \lambda $ such that $\d _\xi\in I'\setminus I$;
thus $\d_\xi \in (\theta ,\theta +\k)$ and, as above, $\cf (\d _\xi )<\k$, which is false because $\,\cf (\d _\xi)=\k$.
If $I'\varsubsetneq I$, then there is $n\in \o$ such that $\d _n\in I\setminus I'$;
thus $\d_n \in (\theta ' ,\theta ' +\k)$ and, hence, $\cf (\d _n )<\k$, which is false since $\,\cf (\d _n)=\k$.
Thus, conditions (B) and (C) are disjoint and in the same way we show that the pairs (B),(D) and (C),(D) are disjoint as well.
\hfill $\Box$
\begin{ex}\label{EX201}\rm
Regarding the classification given in Lemma \ref{T4121} and the ordinals $\d$ of, say, size $\o _2$,
(A) holds for $\o _2 +1$ and $\o _2 +\o$,
(B) for $\o _2 \o +\o _2$,
(C) for $\o _2 \o _1 +\o _2$,
(D) for $\o _2 $ and $\o _2 +\o _2$,
and (E) holds for $\o _2 \o _2 $.
\end{ex}
\paragraph{The separative modification of $\P (\o ^\d)$}
\begin{fac}\label{T4425}
If $\d >0$ is an ordinal, then  we have

(a) $\P (\o ^\d)=\la \I _{\o ^\d }^+ ,\subset\ra$, $\sm (\P (\o ^\d))=\la \I _{\o ^\d }^+ ,\subset_{\I _{\o ^\d} }\ra$ and $\sq (\P (\o ^\d))=(P(\o ^\d )/\I _{\o ^\d })^+$;

(b) If $\d \geq \o$, then $|\sq (\P (\o ^{\d }))|\leq 2^{|\d|}$.
\end{fac}
\dok
(a) By Fact \ref{T4400}(f), $\I _{\o ^\d } =\{ I\subset \o ^\d : \o ^\d \not\hookrightarrow I\}$ is an ideal in $P(\o ^\d )$.
Clearly we have $\P (\o ^\d)\subset \I _{\o ^\d }^+ $. Conversely, if $A\in \I _{\o ^\d }^+$
and $C\in \P(\o ^\d)$, where $C\subset A$,
then  $\o ^\d \hookrightarrow A \hookrightarrow \o ^\d$, which implies that $\A \cong \o ^\d$ and, hence, $A\in \P (\o ^\d)$.
Thus $\P (\o ^\d)= \I _{\o ^\d }^+ $.
The rest follows from Fact \ref{T2226}(d).

(b) By Fact \ref{T4128}(a), $|\o ^{\d }|=|\d|$; so, by (a), $|\sq (\P (\o ^\d ))|\leq |P(\o ^\d )|=2^{|\d|}$.
\hfill $\Box$
\begin{lem}\label{T201}
Let $\d$ be an ordinal. Then

(a) If $\cf (\d )=\k \geq \o$ and $\d =\lim _{\a \rightarrow \k}\d _\a$, then $\cf (\o ^\d )=\k $
and $\o ^\d = \sum _{\a<\k}\o ^{\d _\a}:=\sum _{\a<\k}L_\a $,
where $L_\a\cong \o ^{\d _\a} $, for $\a< \k$, and $L_\a \cap L_\b =\emptyset $, if $\a\neq \b$.

(b) If $\d =\d '+1$, then $\cf (\o ^\d )=\o $ and
$\o ^\d \cong \sum _\o \o ^{\d '}=\sum _{\a<\o}L_\a $,
where $L_\a\cong \o ^{\d '} $, for $\a< \o$, and $L_\a \cap L_\b =\emptyset $, for $\a\neq \b$.
\end{lem}
\dok
(a) By the assumptions $\d$ is a limit ordinal; so, $\o ^\d =\sup \{ \o ^\eta :\eta <\d\}$,
and for $\eta <\d $ there is $\a <\k$ such that $\eta <\d _\a$,
which implies $\o ^\eta <\o ^{\d _\a}<\o ^\d$.
Thus $\la  \o ^{\d _\a} :\a<\k\ra $ is an increasing cofinal sequence in $\o^\d$
and, hence, $\cf (\o ^\d )=\k $.

Let $\theta :=\sum _{\a < \k}\o ^{\d _{\a }}$
and w.l.o.g.\ suppose that $\theta =\sum _{\a < \k} L_\a$, where $L_\a\cong \o ^{\d _\a}$, $\a< \k$, are disjoint.
Clearly we have $\theta \geq  \o ^{\d _{\a }} $, for all $\a < \k$;
so, $\theta \geq  \sup \{ \o ^{\d _{\a }} : \a < \k \}=\o ^\d$.
Suppose that $\o ^\d <\theta$.
Then we would have $\o ^\d\in L_{\a _0}$, for some $\a_0 <\k$;
let $\xi =\min L_{\a _0}$.
By Fact \ref{T4400}(b), $\o ^\d>\xi $ would imply that $\o ^\d \hookrightarrow L_{\a _0}\cong \o^{\d _{\a _0}}$
which is impossible because $\d _{\a _0}<\d$ and, hence, $\o^{\d _{\a _0}}<\o ^\d$.
So $\o ^\d=\xi$
and, hence, $\o ^\d=\sum _{\a< \a _0}\o ^{\d _{\a }}$.
Now, if $\a _0 =\a +1$ we would have $\o ^\d \hookrightarrow L_{\a }\cong \o^{\d _{\a }}$, which is impossible,
and if $\a_0 $ is a limit ordinal we would have $\cf (\o ^\d)\leq \cf (\a _0 )\leq |\a_0 |<\k$,
which is impossible too.
So, $\o ^\d=\theta =\sum _{\a < \k}\o ^{\d _{\a }}$.

(b) If $\d =\d '+1$, then $\o ^\d =\o ^{\d ' +1}=\o ^{\d ' }\o =\sum _\o \o ^{\d '}$
and $\la \o ^{\d '}n:n\in \o\ra$ is an increasing and cofinal sequence in $\o ^\d$, which gives $\cf (\o ^\d )=\o $.
\hfill $\Box$
\begin{fac}\label{T4139}
If $\d>0$ is an ordinal and $\o ^\d\hookrightarrow L=\sum _{i\in I}L_i$,
where $I$ and $L_i$, $i\in I$, are linear orders and $|I| < \cf (\o ^\d )$,
then $\o ^\d\hookrightarrow L_{i }$, for some $i \in I$. 
\end{fac}
\dok
Let $\k=\cf (\o ^\d )$ and let $\la \b _\a : \a < \k  \ra$ be an increasing cofinal sequence in $\o ^\d $.
Then $\o ^\d=\sup \{\b _\a : \a < \k\}$
and for $\a <\k$ we have $g(\b _\a)\in \bigcup _{i\in I}L_i$,
which implies that there is a unique $i _\a\in I$ such that $g(\b _\a)\in L_{i _\a}$.
Clearly, the function $h: \k \rightarrow I$ defined by $h(\a )=i _\a$ is non-decreasing
and, since $|I|<\k=\bigcup _{i\in I}h^{-1}[\{ i\}]$,
by the regularity of $\k$ there is $i ^*\in I$ such that $h^{-1}[\{ i ^* \}]\in [\k ]^{\k }$.
Let $h^{-1}[\{ i ^* \}]=\{ \a _\xi :\xi <\ \k\}$ be an increasing enumeration.
Now, if $\a \in \k \setminus \a _0$ and $\xi \in \k$ where $\a <\a _\xi$,
then $i ^* =h(\a _0)\leq h(\a ) \leq h(\a _\xi) =i ^*$
thus for each $\a \in [\a _0 ,\k )$ we have $h(\a ) =i ^*$,
that is $g(\b _\a)\in L_{i ^*}$.
So, if $\b _{\a _0}\leq \theta <\o ^\d$,
then, by the cofinality of the sequence $\la \b _\a : \a < \k  \ra$ in $\o ^\d$,
there is  $\a\in \k$ such that $\theta \leq \b _\a$,
which implies that $g(\b _{\a _0})\leq g(\theta )\leq g(\b _\a)$.
So, since $g(\b _{\a _0}), g(\b _\a)\in L_{i ^*}$
and $L_{i ^*}$ is a convex suborder of $L$,
we have $g(\theta )\in L_{i ^*}$, for all $\theta \in[\b _{\a _0}, \o ^\d)$.
By Fact \ref{T4400}(b) we have $\o ^\d\cong [\b _{\a _0}, \o ^\d)\cong g[[\b _{\a _0}, \o ^\d)]\subset L_{i ^*}$,
that is $\o ^\d\hookrightarrow L_{i ^*}$.
\hfill $\Box$
\begin{prop}\label{T4103}
Let $\cf (\o ^\d )=\k $ and $\o ^\d=\sum _{\a<\k}\o ^{\d_\a}$, where $\d _\a \leq \d _\b$, whenever $\a <\b <\k$.
If $L:=\sum _{\a<\k}L_\a $, where $L_\a\cong \o ^{\d _\a} $, for $\a< \k$, and $L_\a \cap L_\b =\emptyset $, for $\a\neq \b$, then

(a) Each $C\in \P (L)$ intersects $\k$-many summands $L_\a$;

(b) If $C\subset L$, then $C\in \P (L)$ iff
\begin{equation}\label{EQ4140}\textstyle
\forall \a,\eta < \k \;\; \exists K\in [\k\setminus \eta ]^{<\k} \;\; L_\a \hookrightarrow \bigcup _{\b \in K }L_\b \cap C  ;
\end{equation}

(c) If, in addition, $\cf (\o ^{\d _\a})\geq \k$, for all $\a <\k $,
\begin{equation}\label{EQ217}\textstyle
C\in \P (L) \Leftrightarrow \forall \a < \k \;\; S^\a _C :=\{ \beta <\k :\o ^{\d _\a} \hookrightarrow L_\b \cap C\}\in [\k ]^{\k } .
\end{equation}
\end{prop}
\dok
(a) Let $C\in \P (L)$ and suppose that $C\subset \bigcup _{\a <\eta}L_\a$, for some $\eta < \k$.
By Fact \ref{T4400}(b) we have $L\cong \bigcup _{\a \geq \eta}L_\a$
and, hence, $L+L\hookrightarrow L$,
which is, by Fact \ref{T4102}(a), impossible.

(b) Let $C\in \P (L)$, let $f: L\hookrightarrow L$, where $C=f[L]=\sum _{\a\in \k} f[L_\a]$
and let $\a_0,\eta _0 \in \k$.
By (a) there are $\a _1 <\k$ and $\eta _1\geq \eta_0$ such that $f[L_{\a _1}]\cap L_{\eta _1}\neq \emptyset$
and, hence, $C_1 :=\bigcup _{\a > \a _1}f[L_\a ]\subset \bigcup _{\b \geq \eta _0}L_\b \cap C$.
Since $C_1$ is a final part of $C$ and $C\cong L$,
by Fact \ref{T4400}(b) there is $g: L\hookrightarrow \bigcup _{\b \geq \eta _0}L_\b \cap C$
and, hence, $g[L_{\a _0}]\subset \bigcup _{\b \geq \eta _0}L_\b \cap C$.
If $\b _0\geq \eta _0$, where $g[L_{\a _0 +1}]\cap L_{\b _0}\neq\emptyset$,
then, since $g[L_{\a_0}]<g[L_{\a_0+1}]$, we have
$g[L_{\a_0}]\subset \bigcup _{\eta _0 \leq \b \leq \b_0 }L_\b \cap C$.
Thus $L_{\a_0}\hookrightarrow \bigcup _{\b \in K }L_\b \cap C$,
where $K=[\eta _0, \b _0]\in [\k \setminus \eta _0]^{<\k }$,
and (\ref{EQ4140}) is true.

Assuming that (\ref{EQ4140}) holds,  by recursion we define sequences $\la K_\a :\a\in \k \ra $ and $\la f_\a :\a\in \k \ra $
such that for each $\a ,\a'\in \k$ we have:

(i) $K_\a \in [\k ]^{<\k }$,

(ii) $\a < \a ' \Rightarrow K_\a < K_{\a '} $ and

(iii) $f_\a : L_\a \hookrightarrow \bigcup _{\b \in K_\a}L_\b \cap C$.

\noindent
By (\ref{EQ4140}) for $\a=\eta=0$ there are $K_0 \in [\k ]^{<\k }$ and $f_0 : L_0 \hookrightarrow \bigcup _{\b \in K_0}L_\b \cap C$.

Let $\a '<\k$ and let $K_\a$ and $f_\a$ be defined for $\a <\a '$ such that (i)-(iii) is true.
By (i) and the regularity of $\k $ we have $\eta := \sup (\bigcup _{\a <\a' }K_\a )+1 <\k$
and, by (\ref{EQ4140}), there are $K_{\a '} \in [\k \setminus \eta ]^{<\k }$
and $f_{\a '} : L_{\a '} \hookrightarrow \bigcup _{\b \in K_{\a '}}L_\b \cap C$; the recursion works.

Let $f=\bigcup _{\a \in \k }f_\a$.
By (ii) and (iii), $\a <\a '$ implies $f_{\a }[L_{\a }]<f_{\a '}[L_{\a '}]$
and, hence, $f:L\hookrightarrow C$.
Thus $C'=f[L]\in \P (L)$, $C'\subset C$ and, by Fact \ref{T4425}, $C\in \P (L)$.

(c) By (b) and Fact \ref{T4139}
we have $C\in \P (L)$
iff for each $\a <\k$, for each $\eta < \k$ there is $\beta\geq \eta$ such that $\o ^{\d _\a} \hookrightarrow L_\b \cap C$.
\hfill $\Box$
\paragraph{Complete embeddings of $\CP _\k$}
\begin{te}\label{T205}
If $\d$ is an ordinal, $\cf (\o^\d )=\k\geq \o$,
$\o ^\d  =\sum _{\xi< \k }\o ^{\d _{\xi }}$,
where $\d _\xi \leq \d _\zeta$, whenever $\xi <\zeta <\k$,
and $\cf (\o ^{\d _\xi})\geq\k$, for all $\xi< \k$,
then $\CP _\k \hookrightarrow _c \sq\P (\o ^\d)$.
\end{te}
\dok
Let $\L_\xi\cong \o ^{\d _\xi} $, for $\xi< \k$, where $L_\xi \cap L_\zeta =\emptyset $, for $\xi\neq \zeta$.
Then $\o ^\d \cong \L :=\sum _{\xi < \k} \L_\xi $
and, by Fact \ref{T2226}(d) and Fact \ref{T4425},
$\la P (\o ^\d)\setminus \I _{\o ^\d} ,\subset_{\o ^\d}\ra =\sm \P (\o ^\d)\cong \sm \P (\L)=\la P(L)\setminus \I, \subset _\I\ra$,
where $\I =\{ A\subset L: \L \not\hookrightarrow A \}$.
For $A\subset L$ by Proposition \ref{T4103}(c)  we have
\begin{equation}\label{EQ219}
A\in P(L)\setminus \I \Leftrightarrow \forall \zeta < \k\;\;
S^\zeta_A :=\{ \xi < \k : \o ^{\d _\zeta }\hookrightarrow L_\xi\cap A \}\in [\k ]^\k ,
\end{equation}
and, since $S^\zeta_A\subset \supp (A):=\{\xi< \k : L_\xi\cap  A \neq \emptyset \}$, for $\zeta <\k$,
\begin{equation}\label{EQ220}
|\supp (A) |< \k \Rightarrow A\in \I .
\end{equation}
By Fact \ref{T4432} it is sufficient to  show that $f: \la [\k ]^{\k } ,\subset ^*\ra \hookrightarrow _c \la P(L)\setminus \I , \subset _\I \ra$,
where $f(S)=\bigcup _{\xi\in S}L_\xi$, for $S\in [\k ]^{\k }$.
First, if $S\in [\k ]^{\k }$ and $\zeta < \k$,
then $|S \setminus \zeta|=\k$
and for $\xi \in S \setminus \zeta$ we have $\d _\zeta \leq \d _\xi$,
which gives  $\o ^{\d _\zeta} \hookrightarrow \o ^{\d _\xi}\cong L _\xi=L _\xi \cap f(S)$
and, hence, $\xi \in S^\zeta_{f(S)}$.
Thus $[\k ]^{\k }\ni S \setminus \zeta \subset S^\zeta_{f(S)}$, for each $\zeta < \k$,
and by (\ref{EQ219}) we have $f(S)\in P(L)\setminus \I$.
So,  $f$ maps  $[\k ]^{\k }$ to $P(L)\setminus \I$.
Let $S,T\in [\k ]^{\k }$.

(ce1) If $S\subset ^* T$,
then $|\supp (f(S)\setminus f(T))|=|\supp (\bigcup _{\xi \in S\setminus T }L_\xi )|=|S\setminus T|<\k $
and by (\ref{EQ220}) $f(S)\setminus f(T)\in \I$,
that is $f(S)\subset _\I f(T)$.

(ce2) If $S\perp T$,
then we have $|\supp (f(S)\cap f(T))|=|\supp (\bigcup _{\xi \in S\cap T }L_\xi )|=|S\cap T|<\k$
and by (\ref{EQ220}) $f(S)\cap f(T)\in \I$,
that is $f(S)\perp _{\subset _\I} f(T)$.
If $S\not\perp T$,
then $S\cap T\in [\k ]^{\k } $ and $f(S)\cap f(T)=\bigcup _{\xi\in S\cap T }L_\xi= f(S\cap T)\in P(L)\setminus \I $
and, hence, $f(S)\not\perp _{\subset _\I} f(T)$.

(ce3) First we show that a set $S\in [\k ]^\k$ is a reduction of $A\in P(L)\setminus \I $ to $\la [\k ]^{\k } ,\subset ^*\ra$,
if $S$ is a pseudointersection of the family $\{S^\zeta  _A :\zeta <\k\}$, that is
\begin{equation}\label{EQ4438}
(\forall \zeta < \k \;\; S\subset ^* S^\zeta  _A )\; \Rightarrow \; \forall T\in [\k ]^{\k }\;\; (T \subset ^* S \Rightarrow f(T)\cap  A \in P(L)\setminus \I ) .
\end{equation}
So, let $[\k ]^\k\ni S\subset ^* S^\zeta _A$, for all $\zeta < \k$, and $[\k ]^{\k } \ni T \subset ^* S$.
Then for each $\zeta < \k $ we have $T\subset ^* S^\zeta _A$, which implies that
\begin{equation}\label{EQ221}
\forall \zeta < \k \;\; |T\cap S^\zeta  _A |=\k .
\end{equation}
By recursion we construct a sequence $\la \a_\xi : \xi< \k \ra$ such that for each $\xi < \k $

(i) $\a_\xi \in T$,

(ii) $\xi ' <\xi  \Rightarrow \a_{\xi ' } < \a_{\xi }$ and

(iii) $\o ^{\d _\xi}\hookrightarrow L_{\a_\xi}\cap A$.

\noindent
By (\ref{EQ221}) there is $\a _0\in T \cap S^0_A$ and, hence,  $\o ^{\d _0}\hookrightarrow L_{\a_0}\cap A$.
Let $\zeta <\k $ and let $\la \a_\xi : \xi< \zeta \ra$ be a sequence satisfying (i)--(iii).
By (\ref{EQ221}) we have $|T\cap S^\zeta  _A |=\k$
so there is $\a _\zeta =\min (T\cap S^\zeta  _A \setminus (\sup \{  \a_\xi : \xi< \zeta \}+1))$.
Now we have $\a _\zeta \in T$ and $\a _\zeta>\a _\xi$, for $\xi <\zeta$, so (i) and (ii) are true.
Since $\a _\zeta \in S^\zeta  _A$ we have $\o ^{\d _{\zeta}}\hookrightarrow L_{\a_{\zeta}}\cap A$
so (iii) is true and the recursion works.

Clearly $B:=\bigcup _{\xi< \k } L_{\a_\xi }\cap A\subset A$
and by (i) we have $B\subset \bigcup _{\xi< \k } L_{\a_\xi }\subset \bigcup _{\zeta \in T} L_{\zeta }=f(T)$.
If $\zeta <\k $,
then for $\xi \geq \zeta$ we have $L_{\a_\xi }\cap B=L_{\a_\xi }\cap A$
and, by (iii) and since $\d _\zeta \leq \d _\xi$, we have $\o ^{\d _\zeta }\hookrightarrow \o ^{\d _\xi}\hookrightarrow L_{\a_\xi}\cap A =L_{\a_\xi }\cap B$,
which implies that $\a _\xi \in S^\zeta _B$.
So for each $\zeta <\k$ we have $\{ \a _\xi : \xi \geq \zeta \} \subset S^\zeta _B$,
which by (ii) implies that $S^\zeta _B\in [\k ]^{\k }$.
Thus, by (\ref{EQ219}), $B\in P(L)\setminus \I$
and, since $B\subset f(T)\cap  A  $, we have $f(T)\cap  A \in P(L)\setminus \I$.
So (\ref{EQ4438}) is true.

Finally we prove (ce3). If $A\in P(L)\setminus \I $, then by (\ref{EQ219}) and since $\d _0\leq \d_1 \leq \dots$ we have
\begin{equation}\label{EQ222}
S^0_A \supset S^1_A \supset \dots \supset S^\zeta _A \supset\dots\;\; \mbox{ and }\;\; \forall \zeta <\k \;\; |S^\zeta _A|=\k.
\end{equation}
So, taking $\a _0 =\min S^0_A$
and $\a _\zeta =\min (S^\zeta _A \setminus \{ \a _\xi : \xi <\zeta \})$, for $\zeta <\k$,
we have $S:=\{ \a _\zeta : \zeta <\k \}\in [\k ]^{\k }$
and, by (\ref{EQ222}), $\{ \a _\xi : \xi \geq \zeta \}\subset S^\zeta _A$,
which implies that $S\subset ^* S^\zeta _A$, for all $\zeta <\k $.
By  (\ref{EQ4438}), $S$ is a reduction of $A$ to $\la [\k ]^{\k }, \subset ^* \ra$.
\hfill $\Box$
\begin{te}\label{T4119}
If $\d = \theta + \k$, where $\Ord \ni\theta \geq \k =\cf (\d) >\cf (\theta )=:\lambda \geq \o$,
then $\CP _\lambda \hookrightarrow _c \sq \P (\o ^\d )$.
\end{te}
\dok
By Lemma \ref{T4121}, $\theta =\lim_{\nu \rightarrow \lambda} \d _\nu $, where $\cf (\d _\nu )=\k$, for all $\nu < \lambda$.
Clearly we have $\o ^\d =\o ^{\theta} \o^{ \k}=\sup \{ \o ^{\d _\nu } :\nu < \lambda \}\k$.
Let $\L =\la L, < \ra  =\sum _{\xi < \k } \L_\xi $,
where $\L_\xi=\la L_\xi , <_\xi \ra\cong \o ^{\theta }$, for $\xi < \k$,
and $L_\xi \cap L_\zeta =\emptyset $, for $\xi\neq \zeta$.
Let $\I =\{ A\subset L: \L \not\hookrightarrow A \}$ and, for $A\subset L$ and $\nu < \lambda$, let
$$
S^\nu_A=\{ \xi < \k : \o ^{\d _\nu } \hookrightarrow  L_\xi  \cap A\} .
$$
\begin{cla}\label{T203}
For $A,B\subset L$  we have:

(a) $A\in P(L)\setminus \I \;\Leftrightarrow\;  \forall \nu < \lambda \;\;|S^\nu _A|=\k $;

(b) $A\subset _{\I }B \;\Leftrightarrow\; A\setminus B\in \I\;\Leftrightarrow\; \exists  \nu < \lambda \;\;|S^\nu _{A\setminus B}|<\k $.
\end{cla}
\dok
(a) ($\Rightarrow$) Since $\L \cong \o ^{\theta} \k $, by Proposition \ref{T4103}(b) we have: $A\in P(L)\setminus \I $ iff
\begin{equation}\label{EQ4142'}\textstyle
\forall  \eta  < \k \;\;\exists K \in [\k \setminus \eta  ]^{<\k}\;\;
\o ^\theta \hookrightarrow \bigcup _{\xi\in K }L_\xi \cap A
\end{equation}
and we prove
\begin{equation}\label{EQ4142''}\textstyle
\forall \nu < \lambda \;\;\forall  \eta  < \k \;\;\exists \xi \in \k \setminus\eta  \;\; \xi \in  S^\nu _A.
\end{equation}
If $\nu < \lambda$ and  $\eta  < \k$,
then by (\ref{EQ4142'}) there is $K \in [\k \setminus \eta  ]^{<\k }$
such that $\o ^\theta \hookrightarrow \bigcup _{\xi\in K }L_\xi \cap A$
which implies that $\o ^{\d _\nu } \hookrightarrow \bigcup _{\xi\in K }L_\xi \cap A$
and, since $|K|<\k =\cf(\d _\nu )=\cf(\o ^{\d _\nu })$,
by Fact \ref{T4139} there is $\xi \in K$ such that  $\o ^{\d _\nu } \hookrightarrow L_\xi \cap A$
and, hence, $\xi\in S^\nu _A$.
Since $K \subset \k \setminus \eta $ we have $\xi\geq\eta $.

($\Leftarrow$) Assuming (\ref{EQ4142''}) we prove (\ref{EQ4142'}).
Let $\eta  < \k$.
By recursion we construct a sequence $\la \xi _\nu : \nu < \lambda \ra$
such that for all $\nu, \nu'< \lambda$ we have:

(i) $\nu< \nu' \Rightarrow \eta  \leq \xi _\nu <\xi _{\nu '}<\k$,

(ii) $\o ^{\d _\nu} \hookrightarrow A\cap L_{\xi _\nu }$.

\noindent
First, by (\ref{EQ4142''}), for $\nu =0$ there is $\xi _0 \geq \eta  $ such that $\xi _0 \in  S^0 _A $  that is $\o ^{\d _0} \hookrightarrow A\cap L_{\xi _0}$.
If $\nu < \lambda$ and the sequence $\la \xi _{\nu' } : \nu '< \nu \ra$ satisfies (i) and (ii),
then by (i) and since $\lambda <\k$ we have $\eta  _\nu :=\sup \{\xi _{\nu' } : \nu '< \nu \}+1<\k$
and by (\ref{EQ4142''}) there is $\xi _{\nu } \geq \eta  _\nu$ such that $\xi _{\nu } \in  S^{\nu } _A $
that is $\o ^{\d _{\nu }} \hookrightarrow A\cap L_{\xi _{\nu }}$.
Thus the recursion works.

Now $K:=\{ \xi _\nu :\nu < \lambda \}\in [\k \setminus \eta  ]^{<\k}$
and $\sum _{\nu < \lambda }\o ^{\d _\nu}\hookrightarrow \bigcup _{\nu < \lambda }L_{\xi _\nu } \cap A$,
that is $\o ^\theta \hookrightarrow \bigcup _{\xi\in K}L_{\xi } \cap A$.
So (\ref{EQ4142'}) is true. Statement (b) follows from (a).
\kdok
Let $\{ I_\nu: \nu < \lambda\}\subset [\k ]^\k$ be a partition of $\k$.
For $\nu < \lambda $ and $\xi \in I_\nu$, let $L '_\xi \subset L_\xi$ be the initial part of $\L _\xi$ such that $\L _\xi ' \cong \o ^{\d _\nu}$
and let $C =\bigcup _{\nu < \lambda}\bigcup _{\xi \in I_\nu}L_\xi'$.
Then for each $\nu< \lambda$ the set $S^\nu _C=\{ \xi < \k : \o ^{\d _\nu} \hookrightarrow  L_\xi \cap C \}=\bigcup _{\mu\geq \nu}I_\mu$ is of size $\k$
and, by Claim \ref{T203}(a), $C\in \P (\o ^\d)$.
Thus $\P (C )\cong \P (\o ^\d)$
and by Fact \ref{T4432} it is sufficient to prove that $f:\la [\lambda  ]^{\lambda  }, \subset ^* \ra\hookrightarrow _c \sm \P (C )$,
where for $S\in [\lambda  ]^{\lambda  }$
$$\textstyle
f(S)=C_S:=\bigcup _{\nu\in S}\bigcup _{\xi \in I_\nu}L_\xi'.
$$
First, since $S^\nu_{C_S}=\bigcup _{\mu\in S\setminus \nu}I_\mu$, for $\nu < \lambda$,
by Claim \ref{T203}(a) we have $C_S\in \P (\o ^\d)$
and, hence, $C_S\in \P (C)$ as well.
Let $S,T\in [\lambda  ]^{\lambda  }$.

(ce1) If $S\subset ^* T$,
then, since $f(S)\setminus f(T)=\bigcup _{\nu \in S\setminus T }\bigcup _{\xi \in I_\nu}L_\xi'  $, $|S\setminus T|<\lambda $ and $\lambda$ is a regular cardinal,
for $\nu>\sup (S\setminus T)$ we have $S^\nu_{f(S)\setminus f(T)}=\emptyset$,
which by  Claim \ref{T203}(b) gives $f(S)\subset _\I f(T)$.

(ce2) Since $f(S)\cap f(T)=\bigcup _{\nu  \in S\cap T } \bigcup _{\xi \in I_\nu}L_\xi' $,
by  Claim \ref{T203}(a) we have
$S\perp T$
iff $|S\cap T|<\lambda$
iff $f(S)\cap f(T)\in \I$
iff $f(S)\perp _\I f(T)$.

(ce3) For $A\in \P (C)$ we will find a set $S\in [\lambda  ]^{\lambda  }$
such that $f(T)\cap A \in \P (C)$, whenever $[\lambda  ]^{\lambda  }\ni T\subset ^* S$.
By recursion we construct a sequence $\la \nu_\a:\a < \lambda\ra$ in $\lambda$ such that

(i) $\nu_0=\min \{ \nu \in \lambda : |S^0_A\cap I_\nu |=\k \}$,

(ii) $\nu _{\a+1}=\min \{ \nu \in \lambda \setminus (\nu _\a +1) : |S^{\nu_\a +1}_A\cap I_\nu |=\k \}$, for $\a < \lambda$,

(iii) $\nu _\b=\sup \{ \nu _\a :\a <\b\}$, for $\b \in \Lim \cap \lambda$.

\noindent
Since $A\in \P (\L)$,
by Claim \ref{T203}(a) we have $[\k ]^{\k }\ni S^0_A \subset \bigcup _{\nu < \lambda }I_\nu$
and, since $\k$ is a regular cardinal and $\lambda <\k$, the ordinal $\nu_0$ is well defined by (i).
Let $\b <\lambda$ and let $\la \nu_\a :\a <\b \ra$ be a sequence satisfying (i)--(iii).
If $\b =\a  +1$,
then for $\nu \leq \nu_{\a }$ and $\xi \in I_\nu$ we have $\tp (L_\xi ') \leq \o ^{\d _{\nu_{\a }}}$
and, since $A\subset C$, we have $\o ^{\d _{\nu_{\a } +1}}\not \hookrightarrow L_\xi ' \cap A =L_\xi  \cap A $.
Thus $S^{\nu_{\a } +1}_A \subset \bigcup _{\nu \geq \nu_{\a } +1}I_\nu$
and, as above, there is $\nu \geq \nu_{\a } +1$ such that $|S^{\nu_{\a } +1}_A\cap I_\nu |=\k$.
So $\nu_{\a +1}$ is well defined by (ii).
If $\b $ is a limit ordinal, then $\nu _\b$ is well defined by (iii), because $\lambda$ is a regular cardinal and $\nu_\a <\lambda$, for all $\a <\b$.
Thus, the recursion works and by (ii) and (iii) we have  $S:=\{ \nu_\a:\a< \lambda\}\in [\lambda  ]^{\lambda  }$.

Now, if $T\in [\lambda  ]^{\lambda }$ and $T\subset ^* S$,
then $|T\setminus S|<\lambda$
and, hence, there is $\b< \lambda $ such that $\{ \nu_\a :\a\geq \b \}\subset T$.
So, we have $\bigcup _{\a \geq \b}\bigcup _{\xi \in I_{\nu_\a}}L_\xi' \subset f(T)$,
which gives  $B:=\bigcup _{\a \geq \b}\bigcup _{\xi \in I_{\nu_\a}}L_\xi' \cap A\subset f(T)\cap A$.
For $\a \geq \b$  by (ii) we have $|S^{\nu_\a +1}_A\cap I_{\nu_{\a+1}}|=\k$,
which means that there are $\k$-many $\xi \in I_{\nu_{\a+1}}$ such that $\o ^ {\d_{\nu_\a +1}}\hookrightarrow L_\xi \cap A=L_\xi \cap B$.
Since the sequence $\la \nu_\a +1:\a < \lambda\ra$ is cofinal in $\lambda$,
by Claim \ref{T203}(a) we have $B\in \P (C)$
and, consequently, $f(T)\cap A \in \P (C)$.
\hfill $\Box$
\paragraph{Factorization of factors $\P (\o ^\d)$}
\begin{te}\label{T207}
If $\d>0$ is an ordinal and $\cf (\d )=\k $, then for $\P (\o ^\d)$ we have
\begin{itemize}
\item[\rm (A)] If $\d$ is a successor ordinal or $\cf (\d )=\o$ (that is, $\k \leq \o$), \\
               then $\CP _\o \hookrightarrow _c \sq \P (\o ^\d)$ and $\sq \P (\o ^\d)$ is $\sigma$-closed;
\item[\rm (B)] If $\d = \theta + \k$, where $\Ord \ni\theta\geq \k >\cf (\theta )=\o$\\
               and $\theta =\lim _{n\rightarrow \o}\d _n$, where $\cf (\d _n)=\k$, for all $n\in \o$,\\
               then $\CP _\o \hookrightarrow _c \sq \P (\o ^\d)$ and $\sq \P (\o ^\d )$ is $\sigma$-closed;
\item[\rm (C)] If $\d = \theta + \k$, where $\Ord \ni \theta \geq \k >\cf (\theta )=:\lambda > \o$\\
               and $\theta =\lim _{\xi\rightarrow \,\lambda }\d _\xi $, where $\,\cf (\d _\xi)=\k$, for all $\xi\in \lambda $,\\
               then $\CP _\lambda \hookrightarrow _c \sq \P (\o ^\d )$;
\item[\rm (D)] If $\d =\theta +\k$, where $\Ord \ni \theta\geq \cf (\theta )\geq \k> \o$ or $\theta =0$,\\
               then $\CP _\k \hookrightarrow _c \sq \P (\o ^\d )$;
\item[\rm (E)] If $\d =\lim _{\xi \rightarrow \k} \d _\xi $, where $\cf (\d _\xi)=\k >\o$, for all $\xi< \k$,\\
               then $\CP _\k \hookrightarrow _c \sq \P (\o ^\d)$.
\end{itemize}
In cases (A) and (B) we have $\P (\o ^\d) \equiv_{forc} \CP _\o \ast \pi$, where $\pi$ is a $\CP _\o$-name for an $\o$-distributive forcing.
In other cases forcing $\P (\o ^\d)$ collapses (at least) $\o _2$ to $\o$.
\end{te}
\dok
(A)
Let $\d =\theta + 1$.
For $\theta =0$ we have $\o ^\d=\o$ and $\sq \P (\o )=\CP _\o$.
If $\theta >0$,
then by Lemma \ref{T201}(b) we have $\cf (\o ^\d)=\o$ and $\o ^\d =\sum _\o \o ^{\theta}  $.
Since $\theta \geq 1$ we have $\cf (\o ^\theta)\geq \o $;
so, by Theorem \ref{T205}, $\CP _\o \hookrightarrow _c \sq \P (\o ^\d)$.

If $\cf \d =\o$, let $\d =\sup \{\d _i  : i < \o \}$, where  $\la \d _i  : i < \o \ra$ is an increasing sequence of non-zero ordinals.
Then by Lemma \ref{T201}(a) we have $\cf (\o ^\d)=\o$ and $\o ^\d =\sum _{i<\o}\o ^{\d _i}$.
For $i<\o$ we have $\d _i >0$ and, hence, $\cf (\o ^{\d _i})\geq \o$.
Thus, by Theorem \ref{T205}, $\CP _\o \hookrightarrow _c \sq \P (\o ^\d)$.

Now we prove that the preorder $\sm \P (\o ^\d)$ is $\s$-closed.
By Lemma \ref{T201} we have $\cf (\o^\d )=\o$ and $\o ^\d\cong\sum _{i<\o}\o ^{\d _i}=\sum _{i<\o}L_i=:L$.
Thus $\sm \P (\o ^\d) \cong \sm \P (L)=\la P(L) \setminus \I ,\subset _\I\ra$,
where $\I :=\{ A\subset L :L\not\hookrightarrow A\}$.
For $\d>1$ as in (a) we conclude that $\cf (\o ^{\d _i})\geq \o$, for all $i<\o$,
and, by  Proposition \ref{T4103}(c), for $A\subset L$ we have
\begin{equation}\label{EQ202}\textstyle
A\in P(L) \setminus \I \;\;\Leftrightarrow \;\forall i\in \o \;\; S^i_A:=\{j\in \o : L_i \hookrightarrow L_j \cap A \}\in [\o ]^\o.
\end{equation}
For $A_n \in  P(L) \setminus \I$, $n\in \o$, where $\dots \subset _\I A_1 \subset _\I A_0 $,
we construct $A \in  P(L) \setminus \I$ such that $A\subset _\I A_n$, for all $n\in \o$.
Since $\I$ is an ideal we have $P(L) \setminus \I\ni B_i :=\bigcap _{n\leq i}A_n \subset A_i$, so we can assume that
$A_0 \supset A_1 \supset A_2 \supset \dots$.
By recursion we define
the sequences $\la j_i :i\in \o \ra $ and $\la f_i :i\in \o \ra $ such that for each $i\in \o$ we have

(i) $j_i \in \o $,

(ii) $j_i < j_{i+1} $,

(iii) $f_i : L_i \hookrightarrow L_{j_i} \cap A_i$.

\noindent
Since $A_0\in P(L) \setminus \I$, by (\ref{EQ202}) for $i=0$ there are $j_0 \in \o $ and $f_0 : L_0 \hookrightarrow L_{j_0} \cap C_0$.

Let $j_0, \dots , j_{i}$ and $f_0, \dots , f_{i}$ be sequences satisfying (i)--(iii).
Since $A_{i+1}\in P(L) \setminus \I$
by (\ref{EQ202}) for $i+1$ there is $j_{i+1} \in S^{i+1}_{A_{i+1}} \setminus (j_i+1 )$
so we have $j_i<j_{i+1}$
and there is $f_{i+1} : L_{i+1} \hookrightarrow L_{j_{i+1}} \cap A_{i+1}$.
The recursion works.

By (ii)  $i_1 <i_2$ implies $j_{i_1}<j_{i_2}$,
and, hence, $L_{j_{i_1}}<L_{j_{i_2}}$,
which by (iii) implies $f_{i_1}[L_{i_1}]<f_{i_2}[L_{i_2}]$.
Thus $f=\bigcup _{i\in \o }f_i:L\hookrightarrow L$ and
\begin{equation}\label{EQ4112}\textstyle
A:=f[L]=\bigcup _{i\in \o }f_i[L_i]\in P(L) \setminus \I.
\end{equation}
For $n^* \in \o$ we show that $A\subset _{\I} A_{n^*}$, that is $A\setminus A_{n^*}\in \I$.
By (iii) we have $A\subset \bigcup _{i\in \o}L_{j_i}\cap A_i$
and, since $A_i\subset A_{n^*}$, for $i\geq n^*$,
we have $A\setminus A_{n^*}\subset \bigcup _{i\in \o}L_{j_i}\cap A_i\setminus A_{n^*}=\bigcup _{i<n^*}L_{j_i}\cap A_i\setminus A_{n^*}\subset \bigcup _{i<n^*}L_{j_i}$.
Thus $S^0_{A\setminus A_{n^*}}:=\{ j\in \o : L_0\hookrightarrow L_j \cap A\setminus A_{n^*}\}\subset \{ j_i:i<n^*\}$
and by (\ref{EQ202}) we have $A\setminus A_{n^*}\in \I$.

(B) By Theorem \ref{T4119} (for $\l =\o$) we have $\CP _\o \hookrightarrow _c \sq \P (\o ^\d)$.

By Lemma \ref{T4121}, $\theta =\lim_{n\rightarrow \o} \d _n $, where $\cf (\d _n)=\k$, for all $n\in \o$.
Clearly we have $\o ^\d =\o ^{\theta} \o^{ \k}=\sup \{ \o ^{\d _n} :n\in \o\}\k$.
Let $\L =\la L, < \ra  =\sum _{\xi < \k } \L_\xi $,
where $\L_\xi=\la L_\xi , <_\xi \ra\cong \o ^{\theta }$, for $\xi < \k$,
and $L_\xi \cap L_\zeta =\emptyset $, for $\xi\neq \zeta$.
Let $\I =\{ A\subset L: \L \not\hookrightarrow A \}$ and, for $A\subset L$ and $n\in \o$, let
$$
S^n_A=\{ \xi < \k : \o ^{\d _n} \hookrightarrow  L_\xi \cap A \} .
$$
Then by Claim \ref{T203} (where $\lambda =\o$) for $A,B\subset L$ we have:
\begin{eqnarray}
A\in P(L)\setminus \I & \Leftrightarrow & \forall n\in \o \;\;|S^n _A|=\k ,\label{EQ215}\\
A\subset _{\I }B      & \Leftrightarrow & \exists  n\in \o \;\;|S^n _{A\setminus B}|<\k \label{EQ216}.
\end{eqnarray}
We prove that $\sq \P (\o ^\d )$ is $\sigma$-closed.
Let $A_k \in \P (L)$, for $k\in \o$, where $\dots \subset_\I A_1 \subset_\I A_0$.
W.l.o.g.\ we can suppose that
\begin{equation}\label{EQ4150}
A_0 \supset A_1 \supset \dots ,
\end{equation}
because $\I$ is an ideal
and, hence, $B_k:=\bigcap _{i\leq k}A_k \in \P (L )$ and $B_0 \supset B_1 \supset \dots$;
so a lower bound for $B_k$'s will be a lower bound for $A_k$'s.

Also w.l.o.g.\ we suppose that for each $k\in \o$ we have  $A_k\not\subset _\I A_{k+1}$,
so, by (\ref{EQ216}),  $S^k_{A_k\setminus A_{k+1}}\in [\k ]^{\k }$.
For $\xi \in S^k_{A_k\setminus A_{k+1}}$ we have $\o ^{\d _k} \hookrightarrow L_\xi \cap A_k\setminus A_{k+1}$,
and we choose sets $A^k_\xi$ such that
\begin{equation}\label{EQ4156}
\o ^{\d _k} \cong A^k_\xi \subset  L_\xi \cap A_k\setminus A_{k+1}, \mbox{ for } \xi \in S^k_{A_k\setminus A_{k+1}}\mbox{ and } k\in\o,
\end{equation}
and define
\begin{equation}\label{EQ4157}\textstyle
A=\bigcup _{k\in \o}\bigcup _{\xi \in S^k_{A_k\setminus A_{k+1}}} A^k_\xi .
\end{equation}
By (\ref{EQ4156}) and (\ref{EQ4157}), for $\xi \in S^k_{A_k\setminus A_{k+1}}$ we have $\o ^{\d _k} \hookrightarrow  L_\xi \cap A$,
that is $\xi \in S^k_A$;
so $S^k_{A_k\setminus A_{k+1}}\subset S^k_A$
and, hence, $|S^k_A|=\k $, for all $k\in \o$,
which by (\ref{EQ215}) gives $A\in \P (L)$.

We show that $A\subset _\I A_k$, for $k\in \o$.
By (\ref{EQ4156}) we have $A^k_\xi \subset A_k$, for all $\xi \in S^k_{A_k\setminus A_{k+1}}$
so, by (\ref{EQ4150}) and (\ref{EQ4157}),
\begin{equation}\label{EQ4158}\textstyle
A\setminus A_k \subset \bigcup _{l<k}\bigcup _{\xi \in S^l_{A_l\setminus A_{l+1}}} A^l_\xi .
\end{equation}
Suppose that $\o ^{\d _k} \hookrightarrow L_{\xi _0}\cap  A\setminus A_k $, for some $\xi _0 <\k$.
Then by  (\ref{EQ4158}) and since $A^l_\xi \subset L_\xi$ we have
$\o ^{\d _k} \hookrightarrow\bigcup _{l<k \;\land \;\xi _0\in S^l_{A_l\setminus A_{l+1}}} A^l_{\xi _0}=:B$.
But $B$ is a union of $\leq k$ many summands $A^l_{\xi _0}\cong \o ^{\d _l}$, where $l\leq k-1$
and, hence, $\type (B)\leq \o ^{\d _{k-1}}k <\o ^{\d _{k-1}}\o = \o ^{\d _{k-1}+1}\leq \o ^{\d _k}$
and $\o ^{\d _k} \hookrightarrow B$ is impossible.
Thus for all $\xi  <\k$ we have $\o ^{\d _k} \not\hookrightarrow L_{\xi }\cap A\setminus A_k $,
that is $S^k _{A\setminus A_k}=\emptyset$
and, by (\ref{EQ216}), $A\subset _\I A_k$.

(C) The statement follows from Theorem \ref{T4119} (for $\lambda >\o$).

(D) Clearly we have $\cf (\d )=\k$ so, by Lemma \ref{T201}(a), $\cf (\o^\d )=\k$.
In addition we have $\o ^\d=\o ^{\theta + \k } =\o ^\theta \k =\sum _\k \o^\theta$
and, by Lemma \ref{T201}(a) again, $\cf (\o ^\theta)=\cf (\theta )\geq \k$.
Thus, by Theorem \ref{T205}, $\CP _\k \hookrightarrow _c \sm \P (\o ^\d )$.

(E)
By Lemma \ref{T201}(a) we have $\cf (\o ^\d )=\k$ and $\o ^\d  =\sum _{\xi< \k }\o ^{\d _{\xi }}$.
By the assumption and Lemma \ref{T201}(a) again, for each $\xi< \k$ we have $\cf (\o ^{\d _\xi})=\cf (\d _\xi)\geq\k$;
thus, by Theorem \ref{T205}, $\CP _\k \hookrightarrow _c \sm \P (\o ^\d)$.
\hfill $\Box$
\paragraph{The partial orders $\P (\a )$. ZFC and consistency results}
\begin{te}\label{T226}
If $\a = \o ^{\d _n}s_n+\dots +\o ^{\d _1}s_1+ m$ is an uncountable ordinal presented in the Cantor normal form,
where $n, s_1, \dots , s_n\in \N$,  $\;0< \d_1 < \dots <\d _n$ and $\,m \in \o$,
then, regarding the classification of ordinals given in Lemma \ref{T4121}, we have

(a) If $\d _i$ satisfies (A) or (B), for each $i\leq n$, then the partial order $\sq (\P (\a ))$ is $\s$-closed
and completely embeds ${\CP _\o }^k$, where $k={\sum _{i=1}^n s_i}$;

(b) Otherwise, $\sq \P (\a )$ completely embeds  $\CP _\l $, for some uncountable regular $\l$ and, hence, collapses at least $\o _2$ to $\o$.
If, in addition, $\cc (\CP _\l)=(2^{|\a|})^+$, then $\ro (\sq (\P (\a ) ))\cong \Col (\o ,2^{|\a|})$.
\end{te}
\dok
Let $\la \d'_j :1\leq j\leq k\ra =\la \d_n, \dots , \d _n, \dots , \d_1, \dots , \d _1 \ra$,
where $\d _i$ appears $s_i$ many times, for $i\leq n$.
Then $\a \cong \sum _{1\leq j\leq k}L_j:=L$, where $L_j\cong \o ^{\d'_j}$, for $j\leq k$,
and by Theorem \ref{T200} we have $\sq \P (\a )\cong \prod _{1\leq j\leq k}\sq \P (\o ^{\d'_j})$.

(a) Under the assumptions, by Theorem \ref{T207} for each $j\leq k$
the poset $\sq \P (\o ^{\d'_j})$ is $\sigma$-closed and $\CP_\o \hookrightarrow _c \sq \P (\o ^{\d'_j})$.
By Fact  \ref{T2226}(b) the poset $\sq \P (\a )$ is $\sigma$-closed
and by Fact \ref{T227} ${\CP _\o }^k\hookrightarrow _c \sq \P (\a )$.

(b) Otherwise, by Lemma \ref{T4121} there are a regular $\l \geq \o _1$ and $j\leq k$,
such that $\CP _\l \hookrightarrow _c\sq (\P (\o ^{\d'_j} ))$ and, since  $\sq (\P (\o ^{\d'_j} ))\hookrightarrow _c\sq \P (\a )$
we have $\CP _\l \hookrightarrow _c\sq (\P (\a ))$.
If $\cc (\CP _\l)=(2^{|\a|})^+$,
then by Fact \ref{T209}(e) forcing  by $\sq \P (\a )$ collapses $2^{|\a |}$ to $\o$
and, hence, $|\sq \P (\a )|=2^{|\a |}$.
So, by Fact \ref{T4130},  $\ro (\sq (\P (\a ) ))\cong \Col (\o ,2^{|\a|})$.
\hfill $\Box$
\begin{te} \label{T231}
If $\h <\fc =\o _2= 2^{\o _1}$, then $\ro( P(\o )/\Fin )\cong\Col (\o _1 ,\fc)$ and

(a) For each non-zero ordinal $\d < \fc$ we have
$$
\ro(\sq (\P (\o ^\d ))) \cong \left\{
                     \begin{array}{cl}
                          \Col (\o _1 ,\fc), & \mbox{ if (A) or (B) holds or }\d <\o _1 ,\\[1mm]
                          \Col (\o  ,\fc)                          , & \mbox{ if (D) or (E) holds;}
                     \end{array}
                   \right.
$$

(b) $\ro(\sq (\P (\a ))) \cong \Col (\o _1 ,\fc)$
or $\ro(\sq (\P (\a ))) \cong \Col (\o ,\fc)$, for $\a \in [\o , \fc)$.
\end{te}
\dok
First, by Facts \ref{T209}(b) and \ref{T4130}(a) we have $\ro( P(\o )/\Fin )\cong\Col (\o _1 ,\fc)$.

(a) For $\d <\o _1$ by Theorem \ref{T4421}(b) we have $\ro(\sq (\P (\o ^\d )))\cong \Col (\o _1 ,\fc)$.
If $\d \in [\o _1, \fc )$, then regarding Lemma \ref{T4121} Case (C) impossible.
If $\d$ satisfies (A) or (B), then by Theorem \ref{T207} $\sq \P (\o ^\d )$ is $\o _1$-closed and $\CP _\o \hookrightarrow _c \sq \P (\o ^\d )$,
so, by Fact \ref{T209}(b), $\sq \P (\o ^\d )$ collapses $\fc$ to $\o _1$.
Since $\fc ^{<\o _1}=\fc$, by Fact \ref{T4130} we have  $\ro(\sq (\P (\o ^\d ))) \cong \Col (\o _1 ,\fc)$.
If (D) or (E) holds for $\d$,
then  by Theorem \ref{T207} $\sq \P (\o ^\d )$ collapses $\fc$ to $\o$ and, by Fact \ref{T4130},  $\ro(\sq (\P (\o ^\d ))) \cong \Col (\o  ,\fc)$.

(b) If $\a \in [\o _1 , \fc)$ and $\a = \o ^{\d _n}s_n+\dots +\o ^{\d _1}s_1+ m$,
then by Fact \ref{T4128}(c) we have $|\a|=|\d _n|=\o _1$
and by Fact \ref{T4425}(b) $|\sq (\P (\a ))|=|\sq (\P (\o ^{\d _n} ))|\leq 2^{\o _1}=\fc$.
By Theorem \ref{T226}, either the poset $\sq (\P (\a ))$ is $\s$-closed and collapses $\fc$ to $\o _1$,
or  collapses $\fc$ to $\o$; thus, the statement follows from Fact \ref{T4130}.
\kdok
E.g., regarding Theorem \ref{T231}(a) we obtain $\ro(\sq (\P (\o ^\d ))) \cong\ro( P(\o )/\Fin )$, whenever $\cf (\d )\leq \o$, but also for $\d =\o _1\o +\o _1$.
We note that it is consistent that $\fc =\o _2= 2^{\o _1}$ and $\P (\o ^2 )$ preserves cardinals (see \cite{Kord}).

\section{The posets $\P (\o ^\d)$ and collapsing algebras}\label{S5}
Our intention is to understand how the posets $\P (\a)$ look like.
By Theorem \ref{T200} one step in that direction is
to regard their separative quotients
and to classify the corresponding factors of the form $\sq (\P (\o ^\d))=(P(\o ^\d)/ \I_{\o ^\d})^+$
up to isomorphism of their Boolean completions.
By Theorem \ref{T4421} it is consistent that $\ro(\sq (\P (\o ^\d ))) \cong \Col (\o _1 ,\fc)\cong \ro( P(\o )/\Fin )$
for {\it each} non-zero ordinal $\d < \o _1$ and Theorem \ref{T231}(a) gives an extension for uncountable ordinals.
\paragraph{Collapse to $\o$}
Here we detect some classes of ordinals $\d$ such that $\ro (\sq (\P (\o ^\d )))\cong \Col (\o ,2^{|\d|} )$.
More generally, we have
\begin{fac}\label{T230}
If $\l\geq \o$ is a regular cardinal, $\d \geq \o_1$ an ordinal, $2^{|\d |}>\l$,
and the forcing $\sq(\P (\o ^\d))$ is $\l$-closed and collapses $2^{|\d |}$ to $\l$, then
$\ro (\sq (\P (\o ^\d )))\cong \Col (\l, 2^{|\d |})$.
\end{fac}
\dok
Defining $\B :=P(\o ^\d)/ \I_{\o ^\d}$,
by Fact \ref{T4425} we have $\sq (\P (\o ^\d))=\B ^+$ and $|\B ^+|\leq 2^{|\d|}$.
Since $\B  ^+$ collapses $2^{|\d |}>\l$ to $\l$
we have $|\B ^+|=2^{|\d|}$ and $(2^{|\d|})^{<\l}=2^{|\d|}$;
so, by Fact \ref{T4130}, $\ro (\B )\cong \Col (\l ,2^{|\d |})$.
Clearly, $\sq (\P (\o ^\d))=\B ^+$ implies that $\ro(\sq (\P (\o ^\d)))=\ro(\B )$.
\hfill $\Box$
\begin{te}\label{T222}
$\ro (\sq (\P (\o ^\d )))\cong \Col (\o ,2^{|\d|} )$, if the ordinal $\d$ satisfies (D) or (E),  $2^{\cf (\d )} =2^{|\d|}$ and
($2^{<\cf (\d )}=\cf (\d )$  or $2^{\cf (\d )}=\cf (\d )^+$).
\end{te}
\dok
By Theorem \ref{T207} we have $\k :=\cf (\d ) >\o$ and $\CP _\k \hookrightarrow _c\sq (\P (\o ^\d ))$.
If $2^{<\k}=\k$, then $\cc (\CP _\k)=(2 ^\k)^+$ (since the tree ${}^{<\k}2$ has $2 ^\k$ almost disjoint branches);
by Fact \ref{T209}(a) the equality $2^\k=\k^+$ implies $\cc (\CP _\k)=(2 ^\k)^+$ as well.
Thus, by our assumption, $\cc (\CP _\k)=(2 ^\k)^+=(2^{|\d|})^+$
and by Fact \ref{T209}(e) the poset $\CP _\k ^+$ and, hence, $\sq (\P (\o ^\d ))$, collapses $2^{|\d|}$ to $\o$.
So, by Fact \ref{T230} $\ro (\sq (\P (\o ^\d )))\cong \Col (\o ,2^{|\d|} )$.
\hfill $\Box$
\begin{ex}\rm
We regard Theorem \ref{T222} in the particular case when $|\d |=\cf (\d )=:\k >\o$ and $2^{\cf (\d )} =\cf (\d )^+$.
It is easy to check (using the Cantor normal form (\ref{EQ4146}) for $\d$) that
if (D) holds, then $\d$ satisfies one of the following conditions:
\begin{itemize}\itemsep=-0.7mm
\item[\sc (i)] $\d =\k$
\item[\sc (ii)] $\d= \theta +\k (\g +m)$,\\
      where $\Lim _0 \ni \g <\k $ and $m\geq 2$ and ($\theta =0$ or $\k ^2 \leq \theta <\k ^+$);
\item[\sc (iii)] $\d= \theta +\k ^{\g _2 +m_2} (\g _1 +m _1)+\k$,\\
      where $\Lim _0 \ni \g _1 <\k $ and $\Lim _0 \ni \g _2 <\k ^+$ and $m_1,m_2 \in \N$ and ($\theta =0$ or $\k ^{\g _2 +m_2 +1} \leq \theta <\k ^+$);
\item[\sc (iv)] $\d= \theta +\k ^{\g _2} (\g _1 +m_1 )+\k$,\\
      where $\Lim _0 \ni \g_1 <\k $ and $m_1\in \N$ and  $\cf (\g _2)=\k \leq \g _2 <\k ^+$ and ($\theta =0$ or $\k ^{\g _2 +1} \leq \theta <\k ^+$).
\end{itemize}
On the other hand, if (E) holds, then $\d$ satisfies one of the following conditions:
\begin{itemize}\itemsep=-0.7mm
\item[\sc (i)] $\d= \theta +\k ^{\g +m} $,\\
      where $\Lim _0 \ni \g  <\k ^+$ and $m\in \N$ and $\g +m\leq 2$ and ($\theta =0$ or $\k ^{\g +m} \leq \theta <\k ^+$);
\item[\sc (ii)] $\d= \theta +\k ^\g $,\\
      where $\k =\cf (\g )\leq \g < \k ^+$ and ($\theta =0$ or $\k ^\g \leq \theta <\k ^+$).
\end{itemize}
\end{ex}
If $\d =\k >\o$ is a regular cardinal, then $\sq (\P (\o ^\k ))=\CP _\k$. If $2^{<\k }=\k $ or $2^\k =\k ^+$, then by Theorem \ref{T222} $\ro (\sq (\P (\o ^\k )))\cong \Col (\o ,2^{\k } )$.
\begin{ex}\label{EX203}
It is consistent that  $\ro (\sq (\P (\o ^\k )))\not\cong \Col (\o ,2^{\k } )$, for a regular cardinal $\k >\o$. \rm
Baumgartner \cite{Baum} constructed a model of ZFC in which $\cc (\CP _{\o _1})=\o_3 <2^{\o _1}$.
Thus $\CP _{\o _1 }$ collapses $\o _2$ to $\o$ and the cardinals $> \o _2$ are preserved in extensions by $\CP _{\o _1}$.
So, $\ro (\sq (\P (\o ^{\o_1} )))\not\cong \Col (\o ,2^{\o _1 } )$.
In fact, if $\k >\o$ is a regular cardinal, then  by Fact \ref{T209}(e), $\ro (\sq (\P (\o ^\k )))\cong \Col (\o ,2^\k)$ iff $\cc (\CP _\k)=(2^\k )^+$.
\end{ex}
\begin{te}\label{T228}
$\ro (\sq (\P (\o ^\k )))\cong \Col (\o ,2^\k )$, if $\k$ is a cardinal,  $\k >2^{\cf (\k)}>\cf (\k)>\o$ and $2^\k =\k ^+$.
\end{te}
\dok
Since $\o ^\k=\k$ we have $\sq (\P (\o ^\k ))=\CP _\k^+$
and by Fact \ref{T209}(d) $\sq (\P (\o ^\k ))$ collapses $2^\k =\k ^+$ to $\o$.
So, by Fact \ref{T230}, $\ro (\sq (\P (\o ^\k )))\cong \Col (\o ,2^\k )$.
\hfill $\Box$
\paragraph{Collapse to $\o_1$}
Here we are looking for the ordinals $\d$ such that $\ro (\sq (\P (\o ^\d )))\cong \Col (\o _1,2^{|\d|} )$.
First we show that if $\P (\o ^\d)$ collapses $2^{|\d |}$ to $\o$
(examples of such $\d$'s can be found in the previous paragraph),
then $\ro (\sq (\P (\o ^{\d +n})))\cong \Col (\o _1, 2^{|\d |})$, for all $n\in \N$,
namely, the Boolean completion of $\P (\o ^{\theta})$ is a ``constant function" on the interval $(\d ,\d +\o)$.
We recall some facts from \cite{Kord} and \cite{KCol}.
\begin{fac}\label{T229}
Let $\d$ denote an ordinal, $\k$ a cardinal and $\B\,$ a Boolean algebra.

(a) $\sq (\P (\o ^{\d +n}))\cong (\rp^n (P(\o ^\d)/ \I_{\o ^\d}))^+$, for each  $n\in \N$, whenever $\d \geq 1$.

(b) $\ro (\rp ^n(\Col (\o,\k)))\cong \Col (\o _1, \k )$, for each $n\in \N$, whenever $\k ^\o =\k \geq \fc$.

(c) $\ro (\rp ^n (\Col (\o _1 ,\k)))\cong \Col (\o _1 ,\k ^\o )$, for each $n\in \N$, whenever $\k \geq 2$.

(d) $\ro (\rp ^n(\ro (\B)))\cong \ro (\rp ^n(\B))$, for each $n\in \N$.
\end{fac}
\dok
Statement (a) is a straightforward generalization of Lemma 3.2(d) from \cite{Kord} (proved there for countable $\d$).
Statements (b), (c) and (d) follow from Theorems 6.9 and 5.4 and Fact 2.6 from \cite{KCol}.
\hfill $\Box$
\begin{te}\label{T219}
$\ro (\sq (\P (\o ^{\d +n})))\cong \Col (\o _1, 2^{|\d |})$, for all $n\in \N$, whenever $\d \geq \o_1$ and (a) or (b) holds, where

(a) $\P (\o ^\d)$ collapses $2^{|\d |}$ to $\o$,

(b) $\sq(\P (\o ^\d))$ is $\s$-closed and collapses $2^{|\d |}$ to $\o_1$.
\end{te}
\dok
Defining $\l=\o$, if (a) holds, and $\l=\o_1$, if (b) holds,
by Fact \ref{T230} we have $\ro (\sq (\P (\o ^\d )))\cong \Col (\l , 2^{|\d |})$.
Defining $\B :=P(\o ^\d)/ \I_{\o ^\d}$,
by Fact \ref{T4425} we have $\sq (\P (\o ^\d))=\B ^+$
and, hence, $\ro(\sq (\P (\o ^\d)))=\ro(\B )$;
thus, $\ro(\B )\cong \Col (\l , 2^{|\d |})$.
By Fact \ref{T229}(a) we have $\sq \P (\o ^{\d +n})\cong (\rp ^n (\B))^+$;
so,
\begin{eqnarray*}
\ro(\sq \P (\o ^{\d +n})) & \cong  & \ro(\rp ^n (\B))\quad\mbox{ (since $\rp ^n (\B)$ is a Boolean algebra)}\\
                          & \cong  & \ro(\rp ^n (\ro(\B))) \quad\mbox{ (by Fact \ref{T229}(d))}\\
                          & \cong  & \ro(\rp ^n (\Col (\l ,2^{|\d |}))) \quad\mbox{ (since $\ro(\B )\cong \Col (\l , 2^{|\d |})$)}.
\end{eqnarray*}
Now, if (a) holds, then $\l=\o$ and by Fact \ref{T229}(b) we have $\ro(\rp ^n (\Col (\o ,2^{|\d |})))\cong \Col (\o _1, 2^{|\d |})$.
If (b) holds, then $\l=\o_1$ and, by Fact \ref{T229}(c), $\ro(\rp ^n (\Col (\o _1 ,2^{|\d |})))\cong \Col (\o _1, 2^{|\d |})$.
\hfill $\Box$
\begin{ex}\label{EX202}
It is consistent that $\P (\o ^{\k +1})$ collapses $2^\k$ to $\o _1$, although $\cc (\CP _\k)<(2^\k)^+$ and $\P (\o ^\k)$ does not collapse $2 ^\k$ to $\o $.
\rm
First we note that the equality $\fh =\o _1$ holds in Cohen extensions of models of GCH (see \cite{Blas}).
If $V'$ is a Cohen model for $\fc =\o_5$ obtained from a model $V$ of GCH,
then in $V'$ we have $\fh =\o _1$ and, by Fact \ref{T225}, $2 ^{\o_1}=(\max(\o _5 ,2)^{\o_1})^V=(\o _5^{\o_1})^V=\o _5$.
By Theorem \ref{T207}(A) and Fact \ref{T209}(b) we have $1_{\P (\o ^{\o _1 +1})}\Vdash |2 ^{\o_1}|=\o _1$
although (see \cite{Baum}) $\cc (\CP _{\o _1})=\o _3< \o _6 =(2^{\o _1})^+$ and forcing by $\P (\o ^{\o _1})$ does not collapse $2 ^{\o_1}$ to $\o $.
\end{ex}
\begin{te}\label{T223}
$\ro (\sq (\P (\o ^\mu)))\cong \Col (\o _1 ,2^\mu)$, if $\mu $ is a cardinal, $\mu >\cf (\mu )=\o$ and some of the following conditions hold

(a) $2^\mu =\mu ^+$,

(b) $2^{<\mu}=\mu$,

(c) $\mathrm{MA} + \mu <\fc$,

(d) $V$ is a Cohen model for $\fc =\k$ obtained from a model $V'$ of GCH in which $\k $ is regular and $\mu <\k$.

\noindent
Then $\ro (\sq (\P (\o ^{\mu +n})))\cong \Col (\o _1, 2^\mu)$, for all $n\in \N$.
\end{te}
\dok
By Fact \ref{T4128} we have $\o^\mu =\mu$; thus $\P (\o ^\mu)=\P (\mu )=[\mu ]^\mu$ and by Fact \ref{T2226}(d) we have $\sq \P (\o ^\mu)=\CP _\mu$.
By Theorem \ref{T207} the poset $\CP _\mu$ is $\o _1$-closed
and by Fact  \ref{T209}(c) collapses $\mu ^\o$ to $\o _1$.
So, by Fact \ref{T230} it is sufficient to prove that each of conditions implies that $\mu ^\o =2^\mu$.

If $2^\mu =\mu ^+$,
then, since by K\"{o}nig's theorem $\mu ^+ \leq \mu ^\o \leq2^\mu$, we have $\mu ^\o =2^\mu$.
If  $2^{<\mu}=\mu$, then,
since for a singular cardinal $\mu$ we have $2^\mu \in \{ 2^{<\mu}\mu ^{\cf (\mu )}, (2^{<\mu})^{\cf (2^{<\mu})}\}$,
we obtain $\mu ^\o=2^\mu$.
If $V\models \mathrm{MA} + \mu <\fc $, then $\fc  \leq \mu ^\o \leq 2^\mu=\fc $ and we have $\mu ^\o =2^\mu$.
Finally, if $V=V'_{\Fn (\k ,2)}[G]$,
then by Fact \ref{T225} in $V$ we have
$2^\mu =(\max(\k ,2)^\mu)^{V'}=(\k^\mu)^{V'} =\k$ and
$\mu ^\o=(\max(\k ,\mu)^\o)^{V'}=(\k^\o)^{V'} =\k$;
so, $\mu ^\o =2^\mu$ again.
The second statement follows from the first and Theorem \ref{T219}(b).
\hfill $\Box$
\small
\paragraph{Acknowledgement}
This research was supported by the Science Fund of the Republic of Serbia,
Program IDEAS, Grant No.\ 7750027:
{\it Set-theoretic, model-theoretic and Ramsey-theoretic
phenomena in mathematical structures: similarity and diversity}--SMART.

\end{document}